\documentclass[]{article}

\usepackage{amssymb,latexsym,amsmath,graphicx,subfigure,color,amsthm,IEEEtrantools}     
\usepackage{stmaryrd}
\usepackage{multirow}
\usepackage{epsfig,epstopdf}
\DeclareMathOperator*{\argmin}{argmin}

\newtheorem{theorem}{Theorem}
\newtheorem{lemma}[theorem]{Lemma}

\begin{document}

\title{Explicit and Energy-Conserving Constraint Energy Minimizing Generalized Multiscale Discontinuous Galerkin Method for Wave Propagation in Heterogeneous Media}
\author{
Siu Wun Cheung\thanks{Department of Mathematics, Texas A\&M University, College Station, TX 77843, USA (\texttt{tonycsw2905@math.tamu.edu})}
\and
Eric T. Chung\thanks{Department of Mathematics, The Chinese University of Hong Kong, Shatin, New Territories, Hong Kong SAR, China (\texttt{tschung@math.cuhk.edu.hk})}
\and
Yalchin Efendiev\thanks{Department of Mathematics \& Institute for Scientific Computation (ISC), Texas A\&M University,
College Station, Texas, USA (\texttt{efendiev@math.tamu.edu})}
\and
Wing Tat Leung\thanks{Institute for Computational Engineering and Sciences, 
The University of Texas at Austin, Austin, Texas, USA (\texttt{wleungo@ices.utexas.edu})}
}
\maketitle

\begin{abstract}
In this work, we propose a local multiscale model reduction approach for 
the time-domain scalar wave equation in a heterogenous media. 
A fine mesh is used to capture the heterogeneities of the coefficient field, 
and the equation is solved globally on a coarse mesh in the discontinuous Galerkin discretization setting.  
The main idea of the model reduction approach is to extract dominant modes in local spectral problems 
for representation of important features, 
construct multiscale basis functions in coarse oversampled regions by constraint energy minimization problems, 
and perform a Petrov-Galerkin projection and a symmetrization onto the coarse grid. 
The method is expicit and energy conserving, and 
exhibits both coarse-mesh and spectral convergence, 
provided that the oversampling size is appropriately chosen. 
We study the stability and  convergence of our method. 
We also present numerical results on the Marmousi model in order to test the performance of the method and 
verify the theoretical results.
\end{abstract}

\section{Introduction}
\label{sec:intro}

In this paper, we consider a local multiscale model reduction approach 
for the scalar wave equation. 
Let $\Omega \subset \mathbb{R}^d$ be a bounded spatial domain.
For the sake of simplicity, we focus our discussion on two-dimensional cases in this paper,
though the extension to the three-dimensional case shall be straightforward. 
We consider the following scalar wave equation 
\begin{equation}
\dfrac{\partial^2 u}{\partial t^2} = \text{div} \left(\kappa \nabla u\right) + f \text{ in } [0,T] \times \Omega,
\label{eq:wave}
\end{equation}
where $f(x,t)$ is a given source term. 
The problem is subject to the homogeneous Dirichlet boundary condition $u = 0$ on $[0,T] \times \partial \Omega$, 
and initial conditions $u(x,0) = u_0(x)$ and $u_t(x,0) = v_0(x)$ in $\Omega$.
We assume that the coefficient field $\kappa$ is a heterogeneous 
coefficient field with contrast $\kappa_0 \leq \kappa \leq \kappa_1$. 
Due to the heterogeneities in the coefficient field, 
numerical discretization requires a very fine grid mesh in order to 
capture all the heterogeneities in the medium properties, which potentially
makes 
the numerical solutions on the fine grid become  prohibitively expensive. 
  
Extensive research effort had been devoted to developing numerical solvers 
for solving multiscale problems on the coarse grid, which is typically much coarser than the fine grid,
such as numerical homogenization approaches \cite{papanicolau1978asymptotic,weh02,owhadi2008numerical,abdulle2017multiscale}, 
Multiscale Finite Element Methods (MsFEM) \cite{hw97,ehw99,ch03,eh09,cgh09,jiang2012priori}, 
Variational Multiscale Methods (VMS) \cite{hfmq98,hughes2007variational,Iliev_MMS_11,calo2011note,maier2019explicit,vdovina2005operator}, 
Heterogeneous Multiscale Methods (HMM) \cite{ee03,abdulle05,emz05,engquist2012multiscale,arjmand2014analysis} and 
and Generalized Multsicale Finite Element Methods (GMsFEM) 
\cite{egh12,chung2014adaptive,chung2015generalizedwave,chung2016adaptive,chung2017residualWave,chung2016mixedWave}. 
In numerical homogenization approaches, effective properties are computed for formulating
the global problem on the coarse grid. 
However, these approaches are limited to the cases when the medium properties possess scale separation.
On the other hand, multiscale methods construct of multiscale basis functions which are 
responsible for capturing the local oscillatory effects of the solution. 
Once the multiscale basis functions are constructed, 
coarse-scale equations are formulated.
Moreover, fine-scale information can be recovered by the coarse-scale coefficients and mutliscale basis functions. 
In recent years, multiscale methods in the discontinuous Galerkin (DG) framework have been investigated 
\cite{ehw99,buffa2006analysis,riviere2008discontinuous,eglmsMSDG,elfverson2013dg,efendiev2015spectral,chung2017dg,chung2018dg}. 
In these approaches, unlike conforming finite element formulations, 
multiscale basis functions are in general discontinuous on the coarse grid, 
and stabilization or penalty terms are added to ensure well-posedness of the global problem. 

In many state-of-the-art mutliscale methods, such as MsFEM, VMS and HMM, 
there is 
one basis function per local coarse region to handle the effects of local heterogeneities. 
However, for more complex multiscale problems, each local coarse region contains several high-conductivity regions and 
multiple multiscale basis functions are required to represent the local solution space. 
GMsFEM is developed to allow systematic enrichment of the coarse-scale space with fine-scale information 
and identify the underlying low-dimensional local structures for solution representation. 
The main idea of GMsFEM is to extract local dominant modes by carefully designed local spectral problems in coarse regions, 
and the convergence of the GMsFEM is related to eigenvalue decay of local spectral problems. 
For a more detailed discussion on GMsFEM, we refer the readers to 
\cite{egw10,egh12,eglp13oversampling,chung2014adaptive,chung2015residual,
chung2016adaptive,efendiev2017bayes,cheung2018mc,cheung2019bayes,park2019mc,wang2020vug,maria2019nlmc,cheung2019dg} 
and the references therein. 
Through the design of local spectral problems, our method results in the minimal degree of freedom 
in representing high-contrast features. 
In particular, \cite{chung2014generalized} considered an application of GMsFEM on scalar wave equations. 
On the other hand, several multiscale methods with mesh convergence are developed. 
\cite{owhadi2014polyharmonic, maalqvist2014localization, owhadi2017multigrid}. 
This idea can be adopted for multiscale model reduction techniques for 
achieving both spectral and mesh convergence 
\cite{hou2017sparse,chung2018constraint,chung2018mixed,cheung2018mc,cheung2019dg}.

In this paper, we present the Constraint Energy Minimizing Generalized Multiscale
Discontinuous Galerkin Method (CEM-GMsDGM). 
Our method results in coarse-scale equations 
in an interior penalty discontinuous Galerkin (IPDG) discretization setting.
The method is expicit and energy conserving, and 
exhibits both coarse-mesh convergence and spectral convergence.
The advantages of the method are verified both theoretically and numerically. 
The model reduction approach possesses of two key ingredients.  
The first main ingredient is the local spectral problems in each coarse block for identification of multiscale test basis functions. 
The low-energy dominant modes, which are eigenvectors corresponding to small eigenvalues of local spectral problems, 
are used as multiscale test basis functions, as well as for further construction of 
the second ingredient of our method, which is a set of multiscale trial basis functions. 
Each of the test basis functions sets up an independent orthongonality constraints and 
uniquely defines a corresponding multiscale trial basis function. 
The multiscale trial basis functions will then be  
used for a coarse-scale represenation of the numerical solution. 
We remark that the local spectral problems and the constraint energy minimization problems 
are carefully designed and supported by our analysis.

The paper is organized as follows. In Section~\ref{sec:dg}, we will introduce the notions of grids, 
and essential discretization details such as DG finite element spaces and IPDG formulation on the coarse grid. 
The details of the proposed method, including the construction of basis functions 
and the corresponding systems of linear equations, will be presented in Section~\ref{sec:method}. 
The stability and the convergence of the method will be analyzed in Section~\ref{sec:analysis}.
Numerical results will be provided in Section~\ref{sec:numerical}.
Finally, a conclusion will be given in Section~\ref{sec:conclusion}.

\section{IPDG formulation}
\label{sec:dg}
\begin{figure}[ht!]
\centering
\includegraphics[width=0.5\linewidth]{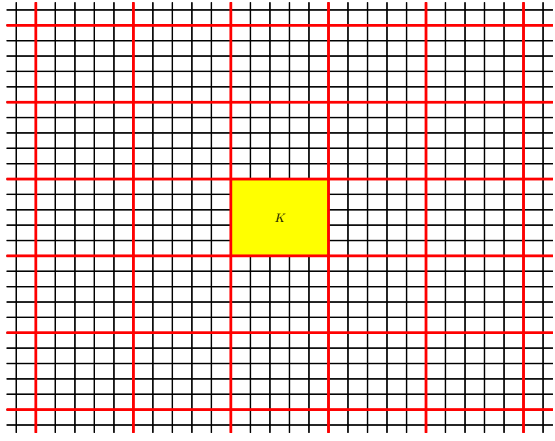}
\caption{An illustration of the fine grid and the coarse grid and a coarse element.}
\label{fig:mesh}
\end{figure}

We are now going to introduce some notions of coarse and fine meshes. 
We start with a usual partition $\mathcal{T}^H$ of $\Omega$ into finite elements, 
which does not necessarily resolve any multiscale features. 
The partition $\mathcal{T}^H$ is called a coarse grid and 
a generic element $K$ in the partition $\mathcal{T}^H$ is called a coarse element. 
Moreover, $H > 0$ is called the coarse mesh size.
We let $N_c$ be the number of coarse grid nodes and 
$N$ be the number of coarse elements. 
We also denote the collection of all coarse grid edges by $\mathcal{E}^H$.
We perform a refinement of $\mathcal{T}^H$ to obtain a fine grid $\mathcal{T}^h$, 
where $h > 0$ is called the fine mesh size. 
It is assumed that the fine grid is sufficiently fine to resolve the solution. 
An illustration of the fine grid and the coarse grid and a coarse element are shown in Figure~\ref{fig:mesh}. 

We are now going to discuss the discontinuous Galerkin (DG) discretization 
and the interior penalty discontinuous Galerkin (IPDG) global formulation. 
For the $i$-th coarse block $K_i$, we 
let $V_h(K_i)$ be 
the conforming bilinear elements defined on the fine grid $\mathcal{T}^h$ in $K_i$.
The DG approximation space is then given by the space of 
coarse-scale locally conforming piecewise bilinear 
fine-grid basis functions, that is, 
\begin{equation}
V_h = \oplus_{i=1}^N V_h(K_i). 
\end{equation}
We remark that functions in $V_h$ are continuous within coarse blocks, 
but discontinuous across the coarse grid edges in general. 
We define the symmetric IPDG bilinear form $a_{{DG}}$ by: 
\begin{equation}
\begin{split}
a_{{DG}}\left(v,w\right) 
& = \sum_{K \in \mathcal{T}^H} \int_K \kappa \nabla v \cdot \nabla w \, dx 
- \sum_{E \in \mathcal{E}^H} \int_E \{ \kappa \nabla v \cdot n_E \} \llbracket w \rrbracket \, d\sigma \\
& \quad - \sum_{E \in \mathcal{E}^H} \int_E \{ \kappa \nabla w \cdot n_E \} \llbracket v \rrbracket \, d\sigma 
+ \dfrac{\gamma}{h} \sum_{E \in \mathcal{E}^H} \int_E \overline{\kappa} \llbracket v \rrbracket \llbracket w \rrbracket \, d\sigma,
\end{split}
\label{eq:dg_bilinear}
\end{equation}
where $\gamma > 0$ is a penalty parameter and 
$n_E$ is a fixed unit normal vector defined on the coarse edge $E \in \mathcal{E}^H$. 
Note that, in \eqref{eq:dg_bilinear}, the average and the jump operators 
are defined in the classical way. 
Specifically, consider an interior coarse edge $E \in \mathcal{E}^H$ and 
let $K^+$ and $K^-$ be the two coarse grid blocks sharing the edge $E$, 
where the unit normal vector $n_E$ is pointing from $K^+$ to $K^-$. 
For a piecewise smooth function $G$ with respect to the coarse grid $\mathcal{T}^H$, we define
\begin{equation}
\begin{split}
\{ G \} & = \dfrac{1}{2}\left(G^+ + G^-\right), \\
\llbracket G \rrbracket & = G^+ - G^-,
\end{split}
\end{equation}
where $G^+ = G\vert_{K^+}$ and $G^- = G\vert_{K^-}$. 
Moreover, on the edge $E$, we define $\overline{\kappa} = \left(\kappa_{K^+} + \kappa_{K^-} \right)/2$, 
where $\kappa_{K^\pm}$ is the maximum value of $\kappa$ over $K^\pm$. 
For a coarse edge $E$ lying on the boundary $\partial \Omega$, we define
$\{G\}=\llbracket G \rrbracket =G$, and $\kappa=\kappa_K$ on $E$, 
where we always assume that $n_E$ is pointing outside of $\Omega$. 
We also use the notation $\left(\cdot,\cdot\right)_{L^2(D)}$ to denote the standard $L^2(D)$ scalar product 
for a subdomain $D \subset \Omega$. 

Using the IPDG spatial discretization, we derive the semi-discrete scheme: 
find $u_h(t,\cdot) \in V_h$ such that 
\begin{equation}
\left(\dfrac{\partial^2 u_h}{\partial t^2}, w\right)_{L^2(\Omega)} + a_{DG}(u_h,w) =  (f,w)_{L^2(\Omega)} \text{ in } [0,T] \times \Omega, 
\label{eq:sol_dg}
\end{equation}
where the initial data is projected onto the finite element space $V_h$ by: 
find $u_h(\cdot,0), \dfrac{\partial u_h}{\partial t}(\cdot,0) \in V_h$ such that for all $w \in V_h$, 
\begin{equation}
\begin{split}
(u_h(\cdot,0), w)_{L^2(\Omega)} & = (u_0, w)_{L^2(\Omega)}, \\
\left(\dfrac{\partial u_h}{\partial t}(\cdot,0), w\right)_{L^2(\Omega)} & = \left(v_0, w\right)_{L^2(\Omega)}. \\
\end{split}
\end{equation} 

\section{Local multiscale model reduction by CEM-GMsFEM}
\label{sec:method}
In this section, we will present our local multiscale model reduction on the IPDG formuation for 
the wave equation by Constraint Energy Minimizing Generlized Multiscale Finite Element Method. 
First, we will use the concept of GMsFEM spectral problems to construct our multiscale test basis functions 
on a generic coarse block $K$ in the coarse grid. 
Next, we will use the concept of constrained energy minimization to construct our multiscale trial basis functions. 
Then, we will derive our coarse-scale model with a Petrov-Galerkin projection and a symmetrc formulation. 
Finally, we present a technique of localization of multiscale trial basis functions on coarse oversampled regions, 
and which results in an explicit time-marching coarse-scale scheme. 

\subsection{Multiscale test functions}
To construct the multiscale test functions, 
we follow the concept of GMsFEM and 
perform a multisale model reduction through a local spectral problem on $V_h(K_i)$, 
which is to find a real number $\lambda_j^{\left(i\right)}$ 
and a function $\phi_j^{\left(i\right)} \in V_h(K_i)$ such that
\begin{equation}
a_i\left(\phi_j^{\left(i\right)}, w\right) = \dfrac{\lambda_j^{\left(i\right)}}{H^2} \left(\phi_j^{\left(i\right)}, w\right)_{L^2(K_i)} 
\text{ for all } w \in V_h(K_i),
\label{eq:spectral_prob}
\end{equation}
where $a_i$ is a symmetric positive semi-definite bilinear form defined as
\begin{equation}
\begin{split}
a_i\left(v,w\right) & = \int_{K_i} \kappa \nabla v \cdot \nabla w \, dx.
\end{split}
\label{eq:spectral_bilinear_form}
\end{equation}
Without loss of generality we shall assume the eigenfunctions are normalized, i.e. 
\begin{equation}
\left(\phi_j^{\left(i\right)}, \phi_{j'}^{\left(i\right)}\right)_{L^2(K_i)} = \delta_{j,j'} \text{ for all } 1 \leq j,j' \leq L_{i}.
\end{equation}
We let $\lambda_j^{\left(i\right)}$ be the eigenvalues of \eqref{eq:spectral_prob} 
arranged in ascending order in $j$, and use the first $L_i$ eigenfunctions 
to construct our local multiscale test space
\begin{equation}
W_H^{\left(i\right)} = \text{span} \{ \phi_j^{\left(i\right)}: 1 \leq j \leq L_i\}.
\end{equation}
We also introduce a local $L^2(K_i)$ projection operator $\pi_i: V_h \to W_H^{\left(i\right)}$ onto $W_H^{\left(i\right)}$ by 
\begin{equation}
\pi_i(v) = \sum_{j=1}^{L_i} \left(v,\phi_j^{\left(i\right)}\right)_{L^2(K_i)} \phi_j^{\left(i\right)} \text{ for all } v \in V_h. 
\end{equation}
The global multiscale test space $W_H$ is then defined as 
the sum of these local multiscale test spaces
\begin{equation}
W_H = \oplus_{i=1}^N W_H^{\left(i\right)}.
\end{equation}
Since the coarse blocks are disjoint, the multiscale test functions 
form an $L^2(\Omega)$ orthonormal basis function for $W_H$, i.e. 
\begin{equation}
\left(\phi_j^{\left(i\right)}, \phi_{j'}^{\left(i'\right)}\right)_{L^2(\Omega)} = \delta_{i,i'} \delta_{j,j'} 
\text{ for all } 1 \leq j \leq L_{i}, 1 \leq j' \leq L_{i'} \text{ and } 1 \leq i,i' \leq N,
\end{equation}
and the global $L^2(\Omega)$ projection operator $\pi: V_h \to W_H$ onto $W_H$ is then naturally defined by 
$\pi = \sum_{i=1}^N \pi_i$.

\subsection{Multiscale trial functions}
Next, we construct our global multiscale trial functions space $V_H$ 
using the concepts of constraint energy minimization.
Given a multiscale test basis function $\psi_{j}^{\left(i\right)}$, 
where $1 \leq j \leq L_i$ and $1 \leq i \leq N$,
the global multiscale trial basis function $\psi_{j}^{\left(i\right)} \in V_h$ is defined as the solution of 
the following constrained energy minimization problem
\begin{equation}
\psi_{j}^{\left(i\right)} = \argmin_{\psi \in V_h} \left\{ a_{{DG}}\left(\psi, \psi\right) : 
\pi(\psi) = \phi_{j}^{\left(i\right)} \right\}.
\label{eq:min1_glo}
\end{equation}
By introducing a Lagrange multiplier, 
the minimization problem \eqref{eq:min1_glo} is equivalent to the following variational problem: 
find $\psi_{j}^{\left(i\right)} \in V_h$ and $\mu_{j}^{\left(i\right)} \in W_H$ such that
\begin{equation}
\begin{split}
a_{{DG}}\left(\psi_{j}^{\left(i\right)}, \psi\right) + \left(\psi, \mu_{j}^{\left(i\right)}\right)_{L^2(\Omega)} & = 0 \text{ for all } \psi \in V_h, \\
\left(\psi_{j}^{\left(i\right)} - \phi_j^{\left(i\right)}, \mu\right)_{L^2(\Omega)} & = 0 \text{ for all } \mu \in W_H.
\end{split}
\label{eq:var1_glo}
\end{equation}
We use the global multiscale trial basis functions to construct
the multiscale trial space, which is defined as
\begin{equation}
\begin{split}
V_{H}^{(\infty)} & = \text{span}  \{\psi_{j}^{\left(i\right)} : 1 \leq j \leq L_i, 1 \leq i \leq N \}.
\label{eq:msdg_space}
\end{split}
\end{equation}

\subsection{Global coarse-scale model}

We derive our fully discrete coarse-scale system by a second-order central difference for temporal discretization.
means of Petrov-Galerkin projection of the fine-scale system onto the coarse-scale spaces.
Let $N_T$ be the number of time steps in the temporal mesh grid and $\tau = T/N_T$ be the time step size. 
At the time instant $t_n = n\tau$, we denote the evaluation of the source function $f$ at the time instant $t_n$ by $f^n$, 
and an approximation of the solution $u(\cdot, t_n)$ by $u_H^n$. 
The coarse-scale model which reads: 
for $n \geq 1$, find $u_H^{n+1} \in V_{H}^{(\infty)}$ such that
\begin{equation}
\left(\dfrac{u_H^{n+1} - 2u_H^n + u_H^{n-1}}{\tau^2}, w \right)_{L^2(\Omega)} + a_{{DG}}\left(u_H^n ,w\right) = \left(f^n, w\right)_{L^2(\Omega)} \,  \text{ for all } w \in W_H, 
\label{eq:sol_ms_pg}
\end{equation}
where the initial data is projected onto the finite element space $V_H^{(\infty)}$ by: 
find $u_h^0, u_h^1 \in V_H^{(\infty)}$ such that for all $w \in V_H^{(\infty)}$, 
\begin{equation}
\begin{split}
(u_H^0, w)_{L^2(\Omega)} & = (u_0, w)_{L^2(\Omega)}, \\
(u_H^1, w)_{L^2(\Omega)} & = \left(u_0 + \tau v_0 + \dfrac{\tau^2}{2} f^0, w\right)_{L^2(\Omega)} - \dfrac{\tau^2}{2} a_{DG}(u_H^0, w).
\end{split}
\end{equation}
Next, we are going to present a symmetric formulation of \eqref{eq:sol_ms_pg} on $V_H^{(\infty)}$. 
By a simple dimensionality argument, 
it is easy to see that $\pi \vert_{V_H^{(\infty)}}: V_H^{(\infty)} \to W_H$ is an isomorphism.
Moreover, for any $v \in V_H^{(\infty)}$, using the fact that 
\begin{equation}
\left( v - \pi(v), w \right)_{L^2(\Omega)} = 0 \text{ for all } w \in W_H, 
\end{equation}
it is straightforward to verify that  
\begin{equation}
a_{DG}(\psi, v - \pi(v)) = 0 \text{ for all } \psi \in V_H^{(\infty)}.
\end{equation}
Combining all these facts, \eqref{eq:sol_ms_pg} can be rewritten as 
\begin{equation}
b\left(\dfrac{u_H^{n+1} - 2u_H^n + u_H^{n-1}}{\tau^2}, w \right) + a_{{DG}}\left(u_H^n ,w\right) = b\left(f^n, w\right) \,  \text{ for all } w \in V_H^{(\infty)},
\label{eq:sol_ms_sym}
\end{equation}
where the bilinear form $b$ is defined as 
\begin{equation}
b(v,w) = \left( \pi(v), \pi(w) \right)_{L^2(\Omega)}. 
\end{equation}

\subsection{Localization of multiscale trial functions}
\begin{figure}[ht!]
\centering
\includegraphics[width=0.5\linewidth]{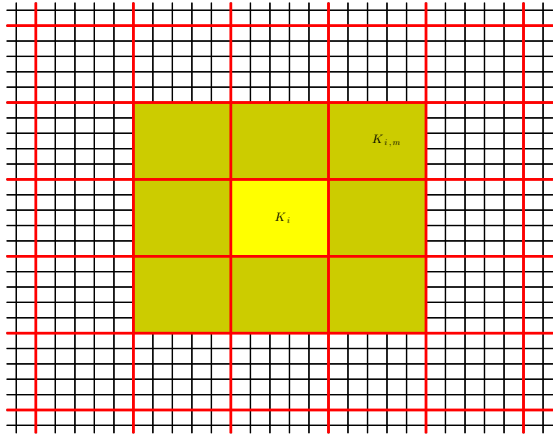}
\caption{An illustration of an oversampled domain formed by enlarging $K_i$ with $1$ coarse grid layer.}
\label{fig:oversample}
\end{figure}

One major drawback of the above constructive procedure 
is that the multiscale trial functions have to be defined 
by solving a global problem. 
Based on our analysis in \cite{cheung2019dg}, 
the global multiscale trial basis function $\psi_{j}^{\left(i\right)}$ exhibits an exponential decay property, 
where the value is very small in locations which are far away from the block $K_i$.
The allows us to construct localized multiscale basis functions on suitably enlarged oversampled domain 
without a significant increase of approximation error. 
More precisely, we denote by $K_{i,m}$ an oversampled domain formed by 
enlarging the coarse grid block $K_i$ by $m$ coarse grid layers, 
of which an illustration shown in Figure~\ref{fig:oversample}.

We introduce the subspaces $V_h\left(K_{i,m}\right)$, which contains restriction of 
fine-scale basis functions in $V_h$ on the oversampled domain $K_{i,m}$, 
and $W_H\left(K_{i,m}\right)$, which is the sum of the local 
multiscale test space over the coarse blocks in the oversampled domain $K_{i,m}$, i.e. 
\begin{equation}
W_H\left(K_{i,m}\right) = \oplus_{K_{i'} \subset K_{i,m}} W_H^{(i')}.
\end{equation}
The localized multiscale trial basis function $\psi_{j,m}^{\left(i\right)} \in V_h$ is then defined as the solution of 
the following constrained energy minimization problem
\begin{equation}
\psi_{j,m}^{\left(i\right)} = \argmin_{\psi \in V_h(K_{i,m})} \left\{ a_{{DG}}\left(\psi, \psi\right) : 
\pi(\psi) = \phi_{j}^{\left(i\right)} \right\}.
\label{eq:min1}
\end{equation}
Using the method of Lagrange multiplier, 
the minimization problem \eqref{eq:min1} is equivalent to the following variational problem: 
find $\psi_{j,m}^{\left(i\right)} \in V_h\left(K_{i,m}\right)$ and $\mu_{j,m}^{\left(i\right)} \in W_H(K_{i,m})$ such that
\begin{equation}
\begin{split}
a_{{DG}}\left(\psi_{j,m}^{\left(i\right)}, \psi\right) + \left(\psi, \mu_{j,m}^{\left(i\right)}\right)_{L^2(K_{i,m})} & = 0 \text{ for all } \psi \in V_h\left(K_{i,m}\right), \\
\left(\psi_{j,m}^{\left(i\right)} - \phi_j^{\left(i\right)}, \mu\right)_{L^2(K_{i,m})} & = 0 \text{ for all } \mu \in W_H(K_{i,m}).
\end{split}
\label{eq:var1}
\end{equation}
The localized multiscale trial basis functions then used to define the localized multiscale trial space, i.e. 
\begin{equation}
\begin{split}
V_H^{(m)} & = \text{span}  \{\psi_{j,m}^{\left(i\right)} : 1 \leq j \leq L_i, 1 \leq i \leq N \}.
\label{eq:msdg_loc_space}
\end{split}
\end{equation}
Finally, our localized coarse-grid model reads: 
for $n \geq 1$, find $u_H^{n+1} \in V_{H}^{(m)}$ such that
\begin{equation}
b\left(\dfrac{u_H^{n+1} - 2u_H^n + u_H^{n-1}}{\tau^2}, w \right) + a_{{DG}}\left(u_H^n ,w\right) = b\left(f^n, w\right) \,  \text{ for all } w \in V_H^{(m)},
\label{eq:sol_ms}
\end{equation}
where the initial data is projected onto the finite element space $V_H^{(m)}$ by: 
find $u_H^0, u_H^1 \in V_H^{(m)}$ such that for all $w \in V_H^{(m)}$, 
\begin{equation}
\begin{split}
(u_H^0, w)_{L^2(\Omega)} & = (u_0, w)_{L^2(\Omega)}, \\
(u_H^1, w)_{L^2(\Omega)} & = \left(u_0 + \tau v_0 + \dfrac{\tau^2}{2} f^0, w\right)_{L^2(\Omega)} - \dfrac{\tau^2}{2} a_{DG}(u_h^0, w). \\
\end{split}
\end{equation}

\subsection{Linear system associated with localized coarse-scale model}
We end this section by the derivation of the resultant coarse-scale system of linear equations 
with a fixed global enumeration of nodal indices and multiscale basis function indices.
Denote by $\mathbf{M}$ and $\mathbf{A}$ the matrix representation of 
the $L^2(\Omega)$ scalar product and the IPDG bilinear form $a_{DG}$ with respect to 
the fine-scale nodal basis functions in $V_h$.
Letting $\mathbf{U}_h(\cdot, t)$ be the column vector consisting of coordinate representation of $u_h(\cdot, t) \in V_h$ 
with respect to fine-scale nodal basis functions in $V_h$, 
the fine-scale system \eqref{eq:sol_dg} can be written as 
\begin{equation}
\mathbf{M} \dfrac{\partial^2 \mathbf{U}_h}{\partial t^2} + \mathbf{A} \mathbf{U}_h = \mathbf{F}.
\end{equation}
It is trivial to see that, for any $v \in V_H^{(m)}$, 
the coefficient representation of $v$ by the multiscale trial basis functions is given by
\begin{equation}
v = \sum_{i=1}^{N} \sum_{j=1}^{L_i} \left(v,\phi_j^{\left(i\right)}\right)_{L^2(\Omega)} \psi_{j,m}^{\left(i\right)},
\end{equation}
which implies the coordinate representations of $v$ with respect to the multiscale trial basis functions $\{ \psi_{j,m}^{(i)} \}$ 
and $\pi(v)$ with respect to the multiscale trial basis functions $\{ \phi_{j}^{(i)} \}$ coincide exactly. 
With a fixed global enumeration of nodal indices and multiscale basis function indices, 
let $\Phi$ and $\Psi$ be the matrix assembled from the column vector consisting the coordinate representation of 
$\phi_j^{(i)}$ and $\psi_{j,m}^{(i)}$ with respect to the fine-scale nodal basis functions in $V_h$. 
Letting $\mathbf{U}_H^{n}$ be the column vector consisting of coordinate representation of $u_H^n \in V_H^{(m)}$ 
with respect to the basis functions $\{ \psi_{j,m}^{(i)} \}$, 
the coarse-scale system \eqref{eq:sol_ms} can be written as 
\begin{equation}
\mathbf{\Phi}^\top \mathbf{M} \mathbf{\Phi} \mathbf{U}_H^{n+1} = 
\mathbf{\Phi}^\top \mathbf{M} \mathbf{\Phi} (2\mathbf{U}_h^{n} - \mathbf{U}_h^{n-1}) +
\tau^2 (\mathbf{\Phi}^\top \mathbf{F}^n - \mathbf{\Psi}^\top\mathbf{A} \mathbf{\Psi} \mathbf{U}_h^n).
\end{equation}
However, noting that the multiscale test basis functions $\{\phi_j^{(i)}\}$ are $L^2(\Omega)$ orthonormal, we have 
$\mathbf{\Phi}^\top \mathbf{M} \mathbf{\Phi} = \mathbf{I}$, and result in an explicit local multiscale model reduction scheme 
\begin{equation}
\mathbf{U}_H^{n+1} = (2\mathbf{U}_H^{n} - \mathbf{U}_H^{n-1}) +
\tau^2 (\mathbf{\Phi}^\top \mathbf{F}^n - \mathbf{\Psi}^\top\mathbf{A} \mathbf{\Psi} \mathbf{U}_h^n).
\end{equation}
Once the explicit scheme is used to obtain the coarse-scale coefficients at a final time step, 
a multiscale approximation is obtained by downscaling with $\mathbf{U}_h^n \approx \mathbf{\Psi} \mathbf{U}_H^{n}$. 

\section{Stability and convergence analysis}
\label{sec:analysis}
In this section, we will analyze the stabiility of proposed localized coarse-grid model \eqref{eq:sol_ms} 
and obtain an error estimate when the coarse-grid solution 
is compared with the fine-grid solution obtained from \eqref{eq:sol_dg}. 
Unless otherwise specified, the constants are generic and 
independent of mesh size and number of oversampling layers. 
First we define a norm on $V_h$ by 
\begin{equation}
\| w \|_{a}^2 = \sum_{K \in \mathcal{T}^H} \int_K \kappa \vert \nabla w \vert^2 \, dx 
+ \dfrac{\gamma}{h} \sum_{E \in \mathcal{E}^H} \int_E \overline{\kappa} \llbracket w \rrbracket^2 \, d\sigma.
\end{equation}
In our analysis, we will make use of following coercivity and continuity results on the IPDG bilinear form, 
provided that the penalty parameter is sufficiently large. 
\begin{lemma}
\cite{chung2014generalized} With a sufficiently large $\gamma > 0$, 
there holds
\begin{equation}
\begin{split}
\vert a_{DG}(v,w) \vert & \leq 2 \| v \|_a \| w \|_a \text{ for all } v, w \in V_h, \\
a_{DG}(v,v) & \geq \dfrac{1}{2} \| v \|_a^2 \text{ for all } v \in V_h. \\
\end{split}
\end{equation}
Moreover, there exists $\beta > 0$ such that 
\begin{equation}
 \| v \|_{L^2(\Omega)}^2 \geq  \beta \kappa_1^{-1} h^2 a_{DG}(v,v) \text{ for all } v \in V_h.
\end{equation}
\label{lemma:ipdg_estimate}
\end{lemma}
With these estimates from the IPDG bilinear form, 
we now establish an inverse Poincar\'{e} inequality on the 
multiscale test function space. 
We will need two fundamental results about from \cite{cheung2019dg}.
The first result is a stability estimate about the projection operator $\pi$.
\begin{lemma}
\cite{cheung2019dg} With a smallness assumption on the fine grid mesh $h>0$, 
there exists $C > 0$ such that 
for any $w_H \in W_{H}$, there exists a function $v \in C^0(\Omega) \cap V_h$ such that 
\begin{equation}
\pi(v) = w_H, \quad 
\| v \|_a^2 \leq D \kappa_1 H^{-2} \| w_H \|_{L^2(\Omega)}^2, \quad
\text{supp}(v) \subseteq \text{supp}(w_H).
\end{equation}
\label{lemma:preimage}
\end{lemma}
The second result states that the global multiscale test functions have a decay property
and their values are small outside a suitably large oversampled domain, 
which in turn justifies the localization of the multiscale test functions. 
\begin{lemma}
\cite{cheung2019dg} With the above assumptions, 
there exists $C, C' > 0$ such that 
\begin{equation}
\| \psi_j^{\left(i\right)} - \psi_{j,m}^{\left(i\right)} \|_a^2 \leq C(1 + C')^{-m}(1+\Lambda^{-1}) \kappa_1 H^{-2} \| \phi_j^{\left(i\right)} \|_{L^2(K_i)}^2, 
\end{equation}
where $\Lambda$ the smallest eigenvalue 
which is excluded in the local spectral problem \eqref{eq:spectral_prob}, i.e. 
\begin{equation}
\Lambda = \min_{1 \leq i \leq N} \lambda_{L_i+1}^{\left(i\right)}.
\end{equation}
\label{lemma:localization}
\end{lemma}
The eigenvalues are related to the contrast in the medium properties. 
In applying to high contrast media, 
the eigenvalues exhibit sharp decay  
and we can simply pick the first few eigenfunctions 
and ensure the smallest eigenvalue excluded is sufficiently large. 
The norm relations in the multiscale space is related to the 
eigenvalue decay in the local spectral problems, 
which accounts for the approximation error by the multiscale finite element space. 
First of all, for any $v \in V_h$, we have 
\begin{equation}
\| (I-\pi)(v) \|_{L^2(\Omega)}^2 \leq 2\Lambda^{-1}H^2 a_{DG}(v,v). 
\label{eq:poincare_ineq}
\end{equation}
Moreover, we have the following inverse inequality. 
\begin{lemma}
With the above assumptions, 
there exists $\beta > 0$ such that 
\begin{equation}
\| \pi(v) \|_{L^2(\Omega)}^2 \geq \beta \kappa_1^{-1} H^2 a_{DG}(v,v) \text{ for all } v \in V_H^{(m)}.
\label{eq:inverse_poincare_ineq}
\end{equation}
\begin{proof}
For any $v \in V_H^{(m)}$, we define
\begin{equation}
\begin{split}
\widetilde{v} & = \sum_{i=1}^N \sum_{j=1}^{L_i} \left(v, \phi_j^{(i)}\right)_{L^2(\Omega)}\psi_{j}^{\left(i\right)}, \\
\eta & = \sum_{i=1}^N \sum_{j=1}^{L_i} \left(v, \phi_j^{(i)}\right)_{L^2(\Omega)} \mu_{j}^{\left(i\right)}.
\end{split}
\end{equation}
Then we have $\pi(v) = \pi(\widetilde{v})$. 
By the definition of global multiscale test functions \eqref{eq:var1_glo}, for any $\psi \in V_h$, we have 
\begin{equation}
\begin{split}
a_{{DG}}\left(\widetilde{v}, \psi\right) + \left(\pi(\psi), \eta\right)_{L^2(\Omega)} & = 0.
\end{split}
\label{eq:var3}
\end{equation}
By Lemma~\ref{lemma:preimage}, we take $w \in C^0(\Omega) \cap V_h$ such that 
\begin{equation}
\pi(w) = \eta, \quad 
\| w \|_a^2 \leq D \kappa_1 H^{-2} \| \eta \|_{L^2(\Omega)}^2, \quad
\text{supp}(w) \subseteq \text{supp}(\eta).
\end{equation}
Taking $\psi = w$ in \eqref{eq:var3}, we have
\begin{equation}
\begin{split}
 \left\|  \eta \right\|_{L^2(\Omega)}^2 
 \leq 2 \left\| \widetilde{v} \right\|_a \| w \|_a 
 \leq 2 D^\frac{1}{2} \kappa_1^\frac{1}{2} H^{-1} \left\| \widetilde{v} \right\|_a \left\|  \eta \right\|_{L^2(\Omega)} 
\end{split}
\end{equation}
Taking $\psi = \widetilde{v}$ in \eqref{eq:var3}, we have
\begin{equation}
\dfrac{1}{2} \left\| \widetilde{v} \right\|_a^2 \leq \left\| \pi(\widetilde{v}) \right\|_{L^2(\Omega)}  \left\|  \eta \right\|_{L^2(\Omega)}.
\end{equation}
Combining these estimates, we have 
\begin{equation}
\left\| \widetilde{v} \right\|_a \leq 4D^\frac{1}{2} \kappa_1^\frac{1}{2} H^{-1} \left\| \pi(v) \right\|_{L^2(\Omega)}.
\end{equation}
On the other hand, by Lemma~\ref{lemma:localization}, we have
\begin{equation}
\begin{split}
\left\| \widetilde{v} - v \right\|_a^2 
& \leq Cm^d \sum_{i=1}^N \left\| \sum_{j=1}^{L_i} \left(v, \phi_j^{(i)}\right)_{L^2(\Omega)} \left(\psi_j^{\left(i\right)} - \psi_{j,m}^{\left(i\right)}\right) \right\|_a^2 \\
& \leq Cm^d (1+\Lambda^{-1}) \kappa_1 H^{-2} (1 + C')^{-m}
\sum_{i=1}^N \left\| \sum_{j=1}^{L_i} \left(v, \phi_j^{(i)}\right)_{L^2(\Omega)} \phi_j^{\left(i\right)} \right\|_{L^2(K_i)}^2 \\
& = C(1+\Lambda^{-1}) m^d (1 + C')^{-m}\kappa_1 H^{-2} \left\| \pi(v) \right\|_{L^2(\Omega)}^2.
\end{split}
\end{equation}
Combining these estimates, we have 
\begin{equation}
\begin{split}
\left\|  v \right\|_a^2 
\leq (16D^2+C(1+\Lambda^{-1}) m^d (1 + C')^{-m}) \kappa_1 H^{-2} \left\| \pi(v) \right\|_{L^2(\Omega)}^2.
\end{split}
\end{equation}
Finally, we note that $m^d (1 + C')^{-m}$ is bounded. The proof is complete by taking 
\begin{equation}
\beta^2 = 16D^2+C (1+\Lambda^{-1})\sup_{m \geq 1} m^d (1 + C')^{-m}.
\end{equation}
\end{proof}
\label{lemma:pi_estimate}
\end{lemma}

We are now going to define a discrete total energy 
which is related to the stability and convergence of our method. 
Given a sequence of states $v = \{v^n\}_{n=0}^{N_T}$, we define the discrete total energy 
at $t = t_{n+\frac{1}{2}}$ by 
\begin{equation}
\begin{split}
E^{n+\frac{1}{2}}(v) & = \dfrac{1}{2} \left\| \pi\left(\dfrac{v^{n+1}-v^n}{\tau}\right) \right\|_{L^2(\Omega)}^2 - 
\dfrac{\tau^2}{8} a_{DG}\left(\dfrac{v^{n+1}-v^n}{\tau},\dfrac{v^{n+1}-v^n}{\tau}\right) \\
&\quad \quad + \dfrac{1}{2} a_{DG}\left(\dfrac{v^{n+1}+v^n}{2},\dfrac{v^{n+1}+v^n}{2}\right),
\end{split}
\end{equation}
which is non-negative under a stability condition. 
More precisely, if there holds 
\begin{equation}
\rho = \dfrac{\kappa_1^{\frac{1}{2}}\tau} {2\beta^{\frac{1}{2}} H} < 1, 
\label{eq:stab_cond}
\end{equation}
then we have the following inequality 
\begin{equation}
E^{n+\frac{1}{2}}(v) \geq \dfrac{1-\rho^2}{2}\left\| \pi\left(\dfrac{v^{n+1}-v^n}{\tau}\right) \right\|_{L^2(\Omega)}^2 + \dfrac{1}{2} a_{DG}\left(\dfrac{v^{n+1}+v^n}{2},\dfrac{v^{n+1}+v^n}{2}\right) \geq 0, 
\label{eq:energy_lb}
\end{equation}
due to the result \eqref{eq:inverse_poincare_ineq} 
in Lemma~\ref{lemma:pi_estimate}. 
The following lemma is the key of proving 
the stability and convergence of our method. 
\begin{lemma}
For $n \geq 1$, given $r^n \in L^2(\Omega)$ 
and $v^{n-1}, v^{n} \in V^{(m)}_H$,
suppose $v^{n+1} \in V^{(m)}_H$ solves
\begin{equation}
b\left(\dfrac{v^{n+1} - 2v^n + v^{n-1}}{\tau^2}, w \right) + a_{{DG}}\left(v^n ,w\right) =
\left(r^n, w\right)_{L^2(\Omega)} \,  \text{ for all } w \in V_H^{(m)}.
\label{eq:sol_v}
\end{equation}
Then we have 
\begin{equation}
E^{n+\frac{1}{2}}(v) = E^{\frac{1}{2}}(v) +\tau  \sum_{k=1}^n \left(r^k, \dfrac{v^{k+1}-v^{k-1}}{2\tau}\right)_{L^2{(\Omega)}}
\label{eq:lem5_eq}
\end{equation}
Moreover, with the above assumptions, then there exist a $C>0$ such that
\begin{equation}
E^{n+\frac{1}{2}}(v) \leq C\left(E^{\frac{1}{2}}(v) +\left(\tau R_1^n +\Lambda^{-\frac{1}{2}}HR_2^n \right)^2\right). 
\label{eq:lem5_ieq}
\end{equation}
where  
\begin{equation}
\begin{split}
R_1^n&=\sum^{n}_{k=1} \|\pi (r^k)\|_{L^2(\Omega)},\\
R_2^n&=\|(I-\pi)r^{1}\|_{L^{2}(\Omega)}+\tau \sum_{k=1}^{n-1}\left\|(I-\pi)\left(\cfrac{r^{k+1}-r^{k}}{\tau}\right)\right\|_{L^{2}(\Omega)}+\|(I-\pi)r^{n}\|_{L^{2}(\Omega)}.
\end{split}
\end{equation}
\begin{proof}
Taking $w = \dfrac{v^{n+1}-v^{n-1}}{2\tau} \in V_H^{(m)}$ in \eqref{eq:sol_v}, we have 
\begin{equation}
\begin{split}
& \dfrac{1}{2\tau}\left(\left\| \pi\left(\dfrac{v^{n+1}-v^n}{\tau}\right) \right\|_{L^2(\Omega)}^2 - 
\left\| \pi\left(\dfrac{v^{n}-v^{n-1}}{\tau}\right) \right\|_{L^2(\Omega)}^2 + 
a_{DG}(v^{n+1}, v^n) - a_{DG}(v^n, v^{n-1})\right)
\\ & \quad 
= \left(r^n, \dfrac{v^{n+1}-v^{n-1}}{2\tau}\right)_{L^2{(\Omega)}}.
\end{split}
\end{equation}
We observe that 
\begin{equation}
\begin{split}
a_{DG}(v^{n+1}, v^n)
& =  a_{DG}\left(\dfrac{v^{n+1}+v^n}{2},\dfrac{v^{n+1}+v^n}{2}\right) - \frac{\tau^2}{4}a_{DG}\left(\dfrac{v^{n+1}-v^n}{\tau},\dfrac{v^{n+1}-v^n}{\tau}\right).
\end{split}
\end{equation}
Hence, we have 
\begin{equation}
\begin{split}
E^{n+\frac{1}{2}}(v) - E^{n-\frac{1}{2}}(v)
= \tau \left(r^n, \dfrac{v^{n+1}-v^{n-1}}{2\tau}\right)_{L^2{(\Omega)}}.
\end{split}
\label{eq:lem5_peq}
\end{equation}
Using a telescoping sum, we obtain \eqref{eq:lem5_eq}.
To obtain the second result, 
we rewrite the right hand side of \eqref{eq:lem5_peq} by
\begin{equation}
\begin{split}
\left(r^{n},\cfrac{v^{n+1}-v^{n-1}}{2\tau}\right)_{L^{2}(\Omega)}
& = b\left(r^{n},\cfrac{v^{n+1}-v^{n-1}}{2\tau}\right)
+\left((I-\pi)r^{n},\cfrac{v^{n+1}-v^{n-1}}{2\tau}\right)_{L^{2}(\Omega)}\\
& = \dfrac{1}{2}\left(b\left(r^{n},\cfrac{v^{n+1}-v^{n}}{\tau}\right)+
b\left(r^{n},\cfrac{v^{n}-v^{n-1}}{\tau}\right)\right)+\\
& \quad \quad \dfrac{1}{\tau}\left((I-\pi)r^{n},(I-\pi)\left(\cfrac{v^{n+1}+v^{n}}{2}\right)\right)_{L^{2}(\Omega)}- \\
& \quad \quad \dfrac{1}{\tau} \left((I-\pi)r^{n},(I-\pi)\left(\cfrac{v^{n}+v^{n-1}}{2}\right)\right)_{L^{2}(\Omega)}.
\end{split}
\label{eq:decompose_rhs}
\end{equation}
Substituting \eqref{eq:decompose_rhs} into \eqref{eq:lem5_eq} and 
rearranging the indices, we obtain 
\begin{equation}
\begin{split}
E^{n+\frac{1}{2}}(v) & = 
E^{\frac{1}{2}}(v) + \dfrac{\tau}{2} \sum_{k=1}^n \left(b\left(r^{k},\cfrac{v^{k+1}-v^{k}}{\tau}\right)+
b\left(r^{k},\cfrac{v^{k}-v^{k-1}}{\tau}\right)\right)+\\
& \quad \quad \left((I-\pi)r^{n},(I-\pi)\left(\cfrac{v^{n+1}+v^{n}}{2}\right)\right)_{L^{2}(\Omega)} - \left((I-\pi)r^{1},(I-\pi)\left(\cfrac{v^{1}+v^{0}}{2}\right)\right)_{L^{2}(\Omega)} - \\
& \quad \quad \tau \sum_{k=1}^{n-1} \left((I-\pi)\left(\dfrac{r^{k+1}-r^k}{\tau}\right),(I-\pi)\left(\cfrac{v^{k+1}+v^{k}}{2}\right)\right)_{L^{2}(\Omega)}. 
\end{split}
\end{equation}
Using Cauchy-Schwarz inequality and Young's inequality, we have 
\begin{equation}
\begin{split}
E^{n+\frac{1}{2}}(v) & \leq 
E^{\frac{1}{2}}(v) +
\max_{0\leq k\leq n} \left\|\pi\left(\cfrac{v^{k+1}-v^{k}}{\tau}\right)\right\|_{L^{2}(\Omega)} \tau R_1^n \\
& \quad \quad + \max_{0\leq k\leq n} \left\|(I-\pi)\left(\cfrac{v^{k+1}+v^{k}}{2}\right)\right\|_{L^{2}(\Omega)}  R_2^n \\
& \leq 
E^{\frac{1}{2}}(v) +
\dfrac{1-\rho^2}{4} \max_{0\leq k\leq n} \left\|\pi\left(\cfrac{v^{k+1}-v^{k}}{\tau}\right)\right\|_{L^{2}(\Omega)}^2 + \dfrac{\tau^2}{1-\rho^2} (R_1^n)^2 \\
& \quad \quad  + \dfrac{\Lambda}{8 H^2} \max_{0\leq k\leq n} \left\|(I-\pi)\left(\cfrac{v^{k+1}+v^{k}}{2}\right)\right\|_{L^{2}(\Omega)}^2+ \dfrac{2 H^2}{\Lambda}(R_2^n)^2.
\end{split}
\end{equation}
Using the inequalities \eqref{eq:poincare_ineq} and \eqref{eq:energy_lb}, 
we obtain the desired result. 
\end{proof}
\label{lemma:energy_estimate}
\end{lemma}

A direct consequence of Lemma~\ref{lemma:energy_estimate} is the following 
stability result of the coarse-grid solution.
\begin{theorem}
With the above assumptions, we have the following stability estimate 
\begin{equation}
\left\|\pi\left(\cfrac{u_{H}^{n+1}-u_{H}^{n}}{\tau}\right)\right\|_{L^{2}(\Omega)}^{2}+
a_{DG}\left(\cfrac{u_{H}^{n+1}+u_{H}^{n}}{2},\cfrac{u_{H}^{n+1}+u_{H}^{n}}{2}\right)\leq 
C\left(E^{\frac{1}{2}}(u_{H})+\tau^2\left(\sum_{k=1}^{n}\|\pi(f^{k})\|_{L^{2}(\Omega)}\right)^{2}\right)
\end{equation}
\end{theorem}

To proceed with our convergence analysis, we need to define two operators. 
The first one is a solution map $G_h: L^2(\Omega) \to V_h$ defined by:
for any $g \in L^2(\Omega)$, the image $G_h g \in V_h$ is defined as 
\begin{equation}
a_{DG}(G_h g, w) = (g,w)_{L^2(\Omega)} \text{ for all } w \in V_h.
\end{equation}
Second, we define an elliptic projection $P_H: V_h \to V_H^{(m)}$ by: 
for any $v \in V_h$, the image $P_H v \in V_H^{(m)}$ is defined as 
\begin{equation}
a_{DG}(P_H v, w) = a_{DG}(v, w) \text{ for all } w \in V_H^{(m)}.
\end{equation}

The approximation error of the DG elliptic projection 
depends on the eigenvalue decay in the local spectral problems 
and has a first-order convergence in the coarse mesh size $H$, 
provided that the discretization follows certain conditions. 
The proof is very much similar to \cite{cheung2019dg} and is omitted. 
\begin{lemma}
\cite{cheung2019dg} With the smallest assumptions 
and the following relation about the size of oversampling region as coarse mesh refines 
\begin{equation}
m = O\left(\log\left(\dfrac{\kappa_1}{H}\right)\right),
\end{equation}
there holds 
\begin{equation}
\| (I-P_H) G_h g \|_a \leq C \Lambda^{-\frac{1}{2}} H \| g \|_{L^2(\Omega)} \text{ for all } g \in L^2(\Omega).
\end{equation}
\label{lemma:a-conv}
\end{lemma}
It is possible to prove the $L^2$ error converges in second order with a duality argument. 
\begin{lemma}
With the above assumptions, there holds 
\begin{equation}
\| (I-P_H) G_h g \|_{L^2(\Omega)} \leq C \Lambda^{-1} H^2 \| g \|_{L^2(\Omega)} \text{ for all } g \in L^2(\Omega).
\end{equation}
\begin{proof}
For any $v \in V_h$, by Galerkin orthogonality, we have 
\begin{equation}
a_{DG}((I-P_H) v, P_H G_h (I-P_H)v) = 0. 
\end{equation}
which implies 
\begin{equation}
\begin{split}
\| (I-P_H) v \|_{L^2(\Omega)}^2 
& = a_{DG}(G_h(I-P_H) v, (I-P_H) v) \\
& = a_{DG}((I-P_H) G_h (I-P_H) v, (I-P_H) v) \\
& \leq 2\| (I-P_H) G_h (I-P_H) v \|_a \| (I-P_H)  v\|_a \\
& \leq C \Lambda^{-\frac{1}{2}} H \| (I-P_H) v \|_{L^2(\Omega)} \| (I-P_H)  v\|_a,
\end{split}
\end{equation}
where we have applied the result from Lemma~\ref{lemma:a-conv} 
with replacing $g$ by $(I-P_H)v$. 
In other words, 
\begin{equation}
\begin{split}
\| (I-P_H) v \|_{L^2(\Omega)}
& \leq C \Lambda^{-\frac{1}{2}} H \| (I-P_H)  v\|_a.
\end{split}
\end{equation}
Applying the result from Lemma~\ref{lemma:a-conv} again with $v = G_h g$, 
we finish our proof. 
\end{proof}
\label{lemma:L2-conv}
\end{lemma}

Now, with all the tools defined, we are going to estimate the error 
between the fine-scale solution $u_h^n = u_h(t_n)$ obtained from solving \eqref{eq:sol_dg}
and the coarse-scale solutions $u_H^n$ obtained from solving \eqref{eq:sol_ms}. 
We define the quantities 
\begin{equation}
\begin{split}
\varepsilon^n & = u_h^n - u_H^n, \\
\delta^n & = u_H^n - P_H u_h^n, \\
\theta^n & = (I-P_H)  u_h^n. \\ 
\end{split}
\end{equation}
We have the following estimates on the 
error of the elliptic projection. 
\begin{lemma}
With the above assumptions and assuming $f \in C^4([0,T]; L^2(\Omega))$, 
there exists $C>0$ such that 
\begin{equation}
\begin{split}
\| \theta^n \|_{L^2(\Omega)} 
& \leq C\Lambda^{-1} H^2 
\left( \left\| g \right\|_{C([0,T];L^2(\Omega))} + 
\tau^2 \left\| \dfrac{\partial^2 g}{\partial t^2} \right\|_{C([0,T];L^2(\Omega))}\right), \\
\left\| \dfrac{\theta^{n+1} - 2\theta^n + \theta^{n-1}}{\tau^2}  \right\|_{L^2(\Omega)} 
& \leq C\Lambda^{-1} H^2 
\left(\left\| \dfrac{\partial^2 g}{\partial t^2} \right\|_{C([0,T];L^2(\Omega))} + 
\tau^2 \left\| \dfrac{\partial^4 g}{\partial t^4} \right\|_{C([0,T];L^2(\Omega))}\right), \\
\| \delta^{1} - \delta^{0}  \|_{L^2(\Omega)} 
& \leq C \tau 
\left(\Lambda^{-1} H^2 \left\| \dfrac{\partial g}{\partial t}  \right\|_{C^0([0,T];L^2(\Omega))} + 
\tau^2 \left\| \dfrac{\partial^3 u_h}{\partial t^3}  \right\|_{C^0([0,T];L^2(\Omega))}\right),
\end{split}
\end{equation}
where $g = f - \dfrac{\partial^2 u_h}{\partial t^2}$.
\label{lemma:diff_err}
\begin{proof}
First, we observe from the definitions that $u_h^n = G_h\left(g(\cdot,t_n) \right)$.
By Taylor's theorem, we have 
\begin{equation}
g(\cdot, t_n + t) = g(\cdot, t_n) 
+ t \dfrac{\partial g}{\partial t}(\cdot, t_n) 
+ \int_{t_n}^{t_n+t} s  \dfrac{\partial^2 g}{\partial t^2}(\cdot, s) ds.
\end{equation}
Integrating from $t = -\tau$ to $t = \tau$, we have 
\begin{equation}
\left\| g(\cdot, t_n) \right\|_{L^2(\Omega)} \leq  
\dfrac{1}{2\tau} \left\| g \right\|_{L^1(t_{n-1},t_{n+1};L^2(\Omega))} + 
\dfrac{\tau}{2} \left\| \dfrac{\partial^2 g}{\partial t^2} \right\|_{L^1(t_{n-1},t_{n+1};L^2(\Omega))}. 
\label{eq:central_error}
\end{equation}
The first result follows directly from Lemma~\ref{lemma:L2-conv}. 
For the second result, 
using Taylor's theorem again, we have 
\begin{equation}
g(\cdot, t_n + t) = g(\cdot, t_n) 
+ t \dfrac{\partial g}{\partial t}(\cdot, t_n) 
+ \dfrac{t^2}{2} \dfrac{\partial^2 g}{\partial t^2}(\cdot, t_n) 
+ \dfrac{t^3}{6} \dfrac{\partial^3 g}{\partial t^3}(\cdot, t_n)
+ \int_{t_n}^{t_n+t} \dfrac{s^3}{6} \dfrac{\partial^4 g}{\partial t^4}(\cdot, s) ds.
\end{equation}
Taking $t = \pm \tau$, we have 
\begin{equation}
\left\| g(\cdot,t_{n+1})-2g(\cdot,t_n)+g(\cdot,t_{n-1}) \right\|_{L^2(\Omega)}
\leq  \tau^2 \left\| \dfrac{\partial^2 g}{\partial t}(\cdot, t_n) \right\|_{L^2(\Omega)} 
+ \dfrac{\tau^3}{3} \left\| \dfrac{\partial^4 g}{\partial t^4} \right\|_{L^1(t_{n-1},t_{n+1};L^2(\Omega))}.
\end{equation}
The second result now follows from Lemma~\ref{lemma:L2-conv} and replacing $g$ by 
$\dfrac{\partial^2 g}{\partial t^2}$ in \eqref{eq:central_error}.  
For the third result, with the same trick, we obtain
\begin{equation}
\left\| g(\cdot,\tau)-g(\cdot,0) \right\|_{L^2(\Omega)}
\leq  \tau \left\| \dfrac{\partial g}{\partial t} \right\|_{C([0,T];L^2(\Omega))},
\end{equation}
which implies 
\begin{equation}
\left \| \theta^1 - \theta^0 \right\|_{L^2(\Omega)} \leq 
C\Lambda^{-1}H^2 \tau \left\| \dfrac{\partial f}{\partial t} -\dfrac{\partial^3 u_h}{\partial t^3}  \right\|_{C^0([0,T];L^2(\Omega))}.
\end{equation}
On the other hand, using Taylor's theorem on $u_h$, we have 
\begin{equation}
u_h(\cdot, \tau) = u_h(\cdot, 0) 
+ \tau \dfrac{\partial u_h}{\partial t}(\cdot, 0) 
+ \dfrac{\tau^2}{2} \dfrac{\partial^2 u_h}{\partial t^2}(\cdot, t_n) 
+ \int_{0}^{\tau} \dfrac{s^2}{2} \dfrac{\partial^3 u_h}{\partial t^3}(\cdot, s) ds.
\end{equation}
Recalling the definition of $u_H^0, u_H^1 \in V_H^{(m)}$ and taking an $L^2(\Omega)$ inner product with $w \in V_H^{(m)}$, we observe that 
\begin{equation}
\begin{split}
\left(u_h^0, w\right)_{L^2(\Omega)} 
& = \left(u_H^0, w\right)_{L^2(\Omega)},  \\
\left(u_h^1, w\right)_{L^2(\Omega)} 
& = \left(u_H^1, w\right)_{L^2(\Omega)}  + \int_{0}^{\tau} \dfrac{s^2}{2} \left(\dfrac{\partial^3 u_h}{\partial t^3}(\cdot, s), w\right)_{L^2(\Omega)} ds.
\end{split}
\end{equation}
This yields 
\begin{equation}
\left \| \varepsilon^1 - \varepsilon^0 \right\|_{L^2(\Omega)} \leq 
\dfrac{\tau^3}{6} \left\| \dfrac{\partial^3 u_h}{\partial t^3}  \right\|_{C^0([0,T];L^2(\Omega))}.
\end{equation}
The third result then follows from a triangle inequality. 
\end{proof}
\label{lemma:theta_err}
\end{lemma}

Now, with all the tools defined, we are going to estimate 
the distance between the coarse-solution solution and elliptic projection. 
\begin{theorem}
With the above assumptions and assuming $f \in C([0,T]; H^1(\Omega))$, 
we have the following estimate
\begin{equation}
\left\|\pi\left(\cfrac{\delta^{n+1}-\delta^{n}}{\tau}\right)\right\|_{L^{2}(\Omega)}^{2}+a_{DG}\left(\cfrac{\delta^{n+1}+\delta^{n}}{2},\cfrac{\delta^{n+1}+\delta^{n}}{2}\right)\leq 
C(\Lambda^{-1}H^{2}+\tau^{2}).
\end{equation}
\begin{proof}
Recalling the definitions \eqref{eq:sol_dg} of $u_h$ and \eqref{eq:sol_ms} of $u_H$ 
and noting that $a_{{DG}}\left(\theta^n ,w\right) = 0$, 
for any $w \in V_H^{(m)}$, we have 
\begin{equation}
b\left(\dfrac{\delta^{n+1} - 2\delta^n + \delta^{n-1}}{\tau^2}, w \right) + a_{{DG}}\left(\delta^n ,w\right) = 
(r^n ,w)_{L^2(\Omega)},
\label{eq:galerkin_diff}
\end{equation}
where $r^n$ is given by
\begin{equation}
r^n = \pi\left(\dfrac{\theta^{n+1} - 2\theta^n + \theta^{n-1}}{\tau^2} \right) + 
\pi\left(\dfrac{\partial^2 u_h}{\partial t^2}(t_{n},\cdot)-\cfrac{u_h^{n+1}-2u_h^{n}+u_h^{n-1}}{\tau^{2}}\right)+(I-\pi)(f^n).
\end{equation}
By Lemma~\ref{lemma:energy_estimate}, we have 
\begin{equation}
\begin{split}
 \left\|\pi\left(\cfrac{\delta^{n+1}-\delta^{n}}{\tau}\right)\right\|_{L^{2}(\Omega)}^{2}+a_{DG}\left(\cfrac{\delta^{n+1}+\delta^{n}}{2},\cfrac{\delta^{n+1}+\delta^{n}}{2}\right)\leq C\left(E^{\frac{1}{2}}(\delta)+\left(\tau R_1^n + \Lambda^{-\frac{1}{2}} HR_2^n\right)^{2}\right),
\end{split}
\end{equation}
where
\begin{equation}
\begin{split}
R_1^n&=\sum_{k=1}^{n}\left\|\pi\left(\cfrac{\theta^{k+1}-2\theta^{k}+\theta^{k-1}}{\tau^2}\right)\right\|_{L^{2}(\Omega)}+\sum_{k=1}^{n}\left\| \pi\left(\dfrac{\partial^2 u_h}{\partial t^2}(t_{k},\cdot)-\cfrac{u_h^{k+1}-2u_h^{k}+u_h^{k-1}}{\tau^{2}}\right) \right\|_{L^{2}(\Omega)},\\
R_2^n&=\|(I-\pi)(f^{1})\|_{L^{2}(\Omega)}+\tau\sum_{k=1}^{n-1}\left\|(I-\pi)\left(\cfrac{f^{k+1}-f^{k}}{\tau}\right)\right\|_{L^{2}(\Omega)}+\|(I-\pi)(f^{n})\|_{L^{2}(\Omega)}.
\end{split}
\end{equation}
First, by Lemma~\ref{lemma:ipdg_estimate} and Lemma~\ref{lemma:pi_estimate}, we note that 
\begin{equation}
E^{\frac{1}{2}}(\delta) \leq \dfrac{1+\rho^2}{2} \left\| \pi\left(\dfrac{\delta^1-\delta^0}{\tau} \right)\right\|_{L^2(\Omega)}^2 + 
 \left\| \dfrac{\delta^1+\delta^0}{2} \right\|_a^2.
\end{equation}
By Lemma~\ref{lemma:a-conv} and the thired inequality in Lemma~\ref{lemma:theta_err}, we have 
\begin{equation}
E^{\frac{1}{2}}(\delta) \leq C(\Lambda^{-1}H^2 + \tau^2). 
\end{equation}
Next, with the second inequality in Lemma~\ref{lemma:theta_err}, 
the first term in $R_1^n$ can be estimated by 
\begin{equation}
\sum_{k=1}^{n}\left\|\pi\left(\cfrac{\theta^{k+1}-2\theta^{k}+\theta^{k-1}}{\tau^2}\right)\right\|_{L^{2}(\Omega)}
\leq C \tau^{-1} \Lambda^{-1} H^2. 
\end{equation}
Similarly, we estimate the second term in $R_1^n$ using Taylor's expansion 
\begin{equation}
\sum_{k=1}^{n} \left\| \pi\left(\dfrac{\partial^2 u_h}{\partial t^2}(t_{k},\cdot)-\cfrac{u_h^{k+1}-2u_h^{k}+u_h^{k-1}}{\tau^{2}}\right) \right\|_{L^{2}(\Omega)}
\leq C \tau \left\| \dfrac{\partial^4 u_h}{\partial t^4} \right\|_{C([0,T]; L^2(\Omega)}.
\end{equation}
On the other hand, we estimate the terms in $R_2^n$ by 
\begin{equation}
\begin{split}
\|(I-\pi)(f^{k})\|_{L^{2}(\Omega)} & \leq \Lambda^{-\frac{1}{2}} H \| f \|_{C([0,T];H^1(\Omega))}, \\
\tau \sum_{k=1}^{n} \left\|(I-\pi)\left(\cfrac{f^{k+1}-f^{k}}{\tau}\right)\right\|_{L^{2}(\Omega)} & \leq C \left\| \dfrac{\partial f}{\partial t} \right\|_{C([0,T];L^2(\Omega)}, 
\end{split}
\end{equation}
Combining all these estimates, we have
\begin{equation}
\tau R_1^n + \Lambda^{-\frac{1}{2}} HR_2^n \leq C (\Lambda^{-1}H^2 + \tau^2), 
\end{equation}
which completes the proof.
\end{proof}
\label{thm:energy_error}
\end{theorem}

We complete this section by providing a $L^2$ error estimate.
\begin{theorem}
With the above assumptions, 
we have the following error estimate
\begin{equation}
\begin{split}
\max_{0 \leq k \leq N_T-1} \left\| \dfrac{\varepsilon^{k+1}+\varepsilon^{k}}{2} \right\|_{L^2(\Omega)}
& \leq C\left(\Lambda^{-1} H^2 + \tau^2 \right).
 \end{split}
\end{equation}
\begin{proof}
We denote $\Delta^n = \tau \sum_{k=1}^n \delta^k$. Using a telescoping sum over \eqref{eq:galerkin_diff}, we have 
\begin{equation}
b\left(\dfrac{\delta^{n+1} - \delta^n }{\tau}, w \right) - b\left(\dfrac{\delta^{1} - \delta^0 }{\tau}, w \right) + a_{{DG}}\left(\Delta^n ,w\right) = 
\tau \sum_{k=1}^n (r^k, w)_{L^2(\Omega)}.  
\end{equation}
Taking $w = \Delta^{n+1} - \Delta^{n-1} = \tau(\delta^{n+1}+\delta^n) \in V_H^{(m)}$, we imply 
\begin{equation}
\begin{split}
& \| \pi (\delta^{n+1}) \|_{L^2(\Omega)}^2 - \| \pi (\delta^{n}) \|_{L^2(\Omega)}^2 + 
 a_{{DG}}\left(\Delta^n ,\Delta^{n+1}\right) - a_{{DG}}\left(\Delta^{n-1} ,\Delta^{n}\right) 
\\ & \quad 
= b\left(\delta^{1} - \delta^0, \delta^{n+1}+\delta^n \right) + 
\tau^2 \sum_{k=1}^n \left(r^k, \delta^{n+1}+\delta^n \right)_{L^2(\Omega)}.
\end{split}
\end{equation}
Using another telescoping sum, we have 
\begin{equation}
\begin{split}
& \| \pi (\delta^{n+1}) \|_{L^2(\Omega)}^2 + 
a_{{DG}}\left(\Delta^n ,\Delta^{n+1}\right) 
- \| \pi (\delta^{1}) \|_{L^2(\Omega)}^2
\\ & \quad 
= \sum_{k=1}^n \left(b\left( \delta^{1} - \delta^0, \delta^{k+1}+\delta^k \right)
+ \tau^2 \sum_{s=1}^k \left(r^s, \delta^{k+1}+\delta^k \right)_{L^2(\Omega)}\right).
\end{split}
\label{eq:error_identity}
\end{equation}
We estimate each of the terms of the error identity \eqref{eq:error_identity}.
For the second term on left hand side of \eqref{eq:error_identity}, we have 
\begin{equation}
\begin{split}
& a_{{DG}}\left(\Delta^n ,\Delta^{n+1}\right) 
\\ & \quad 
= \dfrac{1}{4}\left(a_{{DG}}\left(\Delta^{n+1}+\Delta^{n},\Delta^{n+1}+\Delta^{n}\right)    
- a_{{DG}}\left(\Delta^{n+1}-\Delta^{n}\,\Delta^{n+1}-\Delta^{n}\right)  \right) 
\\ & \quad 
= a_{{DG}}\left(\dfrac{\Delta^{n+1}+\Delta^{n}}{2},\dfrac{\Delta^{n+1}+\Delta^{n}}{2}\right) 
- \dfrac{\tau^2}{4} a_{{DG}}\left(\delta^{n+1},\delta^{n+1}\right) 
\\ & \quad 
 \geq \dfrac{1}{2} \left\| \dfrac{\Delta^{n+1}+\Delta^{n}}{2} \right\|_a^2 - \dfrac{\kappa_1 \tau^2}{4 \beta H^2} \| \pi(\delta^{n+1}) \|_{L^2(\Omega)}^2
 \\ & \quad 
 = \dfrac{1}{2} \left\| \dfrac{\Delta^{n+1}+\Delta^{n}}{2} \right\|_a^2 -\rho^2 \| \pi(\delta^{n+1}) \|_{L^2(\Omega)}^2
\end{split}
\end{equation} 
For the last term on left hand side of \eqref{eq:error_identity}, we proceed with the standard procedure 
with Cauchy-Schwarz inequality to see that 
\begin{equation}
\begin{split}
\| \pi(\delta^1) \|_{L^2(\Omega)}^2 
& = \| \pi(\delta^0) \|_{L^2(\Omega)}^2 + b(\delta^{1} - \delta^{0}, \delta^{1}) + b(\delta^{1} - \delta^{0},\delta^{0}) \\
& \leq \| \pi(\delta^0) \|_{L^2(\Omega)}^2 + \| \pi(\delta^{1} - \delta^{0}) \|_{L^2(\Omega)} 
(\| \pi(\delta^{1}) \|_{L^2(\Omega)} + \| \pi(\delta^{0}) \|_{L^2(\Omega)}) \\ 
& \leq \| \delta^0 \|_{L^2(\Omega)}^2 + 2\| \delta^{1} - \delta^{0} \|_{L^2(\Omega)} 
\max_{0 \leq k \leq N_T} \| \pi(\delta^{k}) \|_{L^2(\Omega)}.
\end{split}
\end{equation}
Similarly, for the first term on the right hand side of \eqref{eq:error_identity}, we have 
\begin{equation}
\begin{split}
b\left(\delta^{1} - \delta^{0}, \delta^{k+1}+\delta^k \right) & 
\leq 2 \| \delta^{1} - \delta^{0} \|_{L^2(\Omega)} \max_{0 \leq k \leq N_T} \| \pi(\delta^{k}) \|_{L^2(\Omega)}. \\
\end{split}
\end{equation}
Finally, we estimate the second term on the right hand side of \eqref{eq:error_identity} by 
\begin{equation}
\begin{split}
& \sum_{k=1}^n \sum_{s=1}^k \left(r^s, \delta^{k+1}+\delta^k\right)_{L^2(\Omega)} \\
& \quad \leq \sum_{k=1}^n R_1^k \left\| \pi(\delta^{k+1}+\delta^k) \right\|_{L^2(\Omega)} + 
 \sum_{k=1}^n R_2^k \left\| (I- \pi)(\delta^{k+1}+\delta^k) \right\|_{L^2(\Omega)}  \\
 & \quad \leq 2 \left(\sum_{k=1}^n R_1^k\right) \max_{0 \leq k \leq N_T} \left\| \pi(\delta^k) \right\|_{L^2(\Omega)} + 
2 \left( \sum_{k=1}^n R_2^k\right) 
\max_{0 \leq k \leq N_T-1} \left\| (I- \pi)\left(\dfrac{\delta^{k+1}+\delta^k}{2}\right) \right\|_{L^2(\Omega)}.
\end{split}
\end{equation}
Combining these estimates, we infer from \eqref{eq:error_identity} that 
\begin{equation}
\begin{split}
& (1-\rho^2) \| \pi(\delta^{n+1}) \|_{L^2(\Omega)}^2 + 
 \dfrac{1}{2} \left\| \dfrac{\Delta^{n+1}+\Delta^{n}}{2} \right\|_a^2
 \\ & \quad 
\leq \| \delta^{0} \|_{L^2(\Omega)}^2 + 2\left((n+1) \| \delta^{1}-\delta^0 \|_{L^2(\Omega)} + \tau^2 \sum_{k=1}^{n} R_1^k \right) 
\max_{0 \leq k \leq N_T} \| \pi(\delta^{k}) \|_{L^2(\Omega)} \\
& \quad \quad \quad + 2 \left(\tau^2 \sum_{k=1}^n R_2^k \right) 
\max_{0 \leq k \leq N_T-1}  \left\| (I- \pi)\left(\dfrac{\delta^{k+1}+\delta^k}{2}\right) \right\|_{L^2(\Omega)}.
\end{split}
\end{equation}
Using Young's inequality, we have 
\begin{equation}
\begin{split}
& \max_{0 \leq k \leq N_T} \| \pi(\delta^{k}) \|_{L^2(\Omega)}^2 
 \\ & \quad 
\leq \dfrac{2}{1-\rho^2} \| \delta^{0} \|_{L^2(\Omega)}^2 + \dfrac{8}{(1-\rho^2)^2} \left(N_T \| \delta^{1}-\delta^0 \|_{L^2(\Omega)} + \tau^2 \sum_{k=1}^{N_T-1} R_1^k \right)^2 \\
& \quad \quad \quad + \dfrac{4}{(1-\rho^2)^2} \left(\tau^2  \sum_{k=1}^{N_T-1} R_2^k\right)^2 + 
\max_{0 \leq k \leq N_T-1} \left\| (I- \pi)\left(\dfrac{\delta^{k+1}+\delta^k}{2}\right) \right\|_{L^2(\Omega)}^2,
\end{split}
\end{equation}
which further implies 
\begin{equation}
\begin{split}
& \max_{0 \leq k \leq N_T-1} \left\| \dfrac{\delta^{k+1}+\delta^{k}}{2} \right\|_{L^2(\Omega)}^2 
 \\ & \quad 
\leq \dfrac{2}{1-\rho^2} \| \delta^{0} \|_{L^2(\Omega)}^2 + \dfrac{8}{(1-\rho^2)^2} \left( N_T \| \delta^{1}-\delta^0 \|_{L^2(\Omega)} + \tau^2 \sum_{k=1}^{N_T-1} R_1^k \right)^2 \\
& \quad \quad \quad + \dfrac{4}{(1-\rho^2)^2} \left(\tau^2 \sum_{k=1}^{N_T-1} R_2^k\right)^2 + 
2 \max_{0 \leq k \leq N_T-1} \left\| (I- \pi)\left(\dfrac{\delta^{k+1}+\delta^k}{2}\right) \right\|_{L^2(\Omega)}^2,
\end{split}
\end{equation}
From the third inequality of Lemma~\ref{lemma:theta_err}, we have 
\begin{equation}
N_T \| \delta^{1}-\delta^0 \|_{L^2(\Omega)} \leq C(\Lambda^{-1}H^2 + \tau^2).
\end{equation}
Applying the estimates on $R_1^k$ and $R_2^k$ in Theorem~\ref{thm:energy_error}, we have 
\begin{equation}
\tau^2 \sum_{k=1}^{N_T-1} (R_1^k + R_2^k) \leq C\left(\Lambda^{-1}H^2 + \tau^2\right).
\end{equation}
Using \eqref{eq:poincare_ineq} and the result of Theorem~\ref{thm:energy_error}, we have 
\begin{equation}
\max_{0 \leq k \leq N_T-1} \left\| (I- \pi)\left(\dfrac{\delta^{k+1}+\delta^k}{2}\right) \right\|_{L^2(\Omega)} 
\leq C (\Lambda^{-1}H^2 + \tau^2). 
\end{equation}
Combining these estimates, we have 
\begin{equation}
\max_{0 \leq k \leq N_T-1} \left\| \dfrac{\delta^{k+1}+\delta^{k}}{2} \right\|_{L^2(\Omega)} 
\leq C(\| \delta^{0} \|_{L^2(\Omega)} + \Lambda^{-1}H^2 + \tau^2).
\end{equation}
Using a triangle inequality with $\delta^0 = \theta^0 - \varepsilon^0$, we have 
\begin{equation}
\max_{0 \leq k \leq N_T-1} \left\| \dfrac{\delta^{k+1}+\delta^{k}}{2} \right\|_{L^2(\Omega)} 
\leq C\left(\| \varepsilon^{0} \|_{L^2(\Omega)} + \max_{0 \leq k \leq N_T} \left\| \theta^k \right\|_{L^2(\Omega)} + \Lambda^{-1}H^2 + \tau^2\right).
\end{equation}
Using another triangle inequality with $\varepsilon^k = \theta^k - \delta^k$, we have 
\begin{equation}
\begin{split}
\max_{0 \leq k \leq N_T-1} \left\| \dfrac{\varepsilon^{k+1}+\varepsilon^{k}}{2} \right\|_{L^2(\Omega)} 
& \leq \max_{0 \leq k \leq N_T-1} \left\| \dfrac{\delta^{k+1}+\delta^{k}}{2} \right\|_{L^2(\Omega)} + \max_{0 \leq k \leq N_T} \left\| \theta^k \right\|_{L^2(\Omega)} \\
& \leq C\left(\| \varepsilon^{0} \|_{L^2(\Omega)} + \max_{0 \leq k \leq N_T} \left\| \theta^k \right\|_{L^2(\Omega)} + \Lambda^{-1}H^2 + \tau^2\right).
\end{split}
\end{equation}
Since $u_H^{0}$ is the $L^2$ projection of $u_h^0$ onto $V_H^{(m)}$, we have 
\begin{equation}
\| \varepsilon^{0} \|_{L^2(\Omega)} \leq \| \theta^{0} \|_{L^2(\Omega)} \leq C\Lambda^{-1}H^2,
\end{equation}
thanks to Lemma~\ref{lemma:L2-conv}. 
The proof is completed by applying the first inequality of Lemma~\ref{lemma:theta_err}.
\end{proof}
\end{theorem}

\section{Numerical results}
\label{sec:numerical}

\begin{figure}[ht!]
\centering
\includegraphics[width=0.6\linewidth]{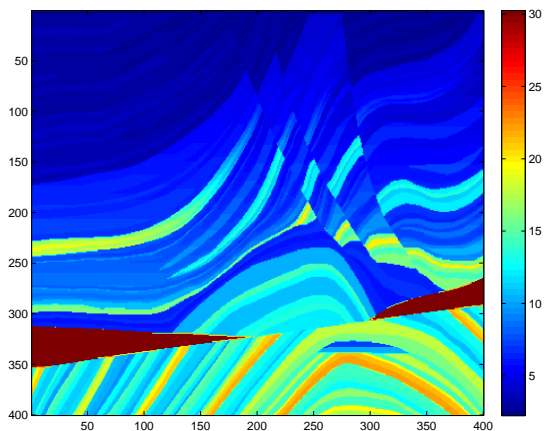}
\caption{Marmousi model for Experiment 2.}
\label{fig:marmousi}
\end{figure}

In this section, we will present numerical examples 
on the scalar wave equation to demonstrate the convergence of our proposed method 
with respect to the coarse mesh size $H$ and 
the number of oversampling layers $m$. 
We take the bulk modulus on the spatial domain $\Omega = [0,1]^2$ 
as part of the Marmousi model as shown in Figure~\ref{fig:marmousi}. 
In all the experiments, the IPDG penalty parameter in \eqref{eq:dg_bilinear} 
is set to be $\gamma = 4$, 
which is experimentally sufficient for ensuring the coercivity of the bilinear form $a_{DG}$.
The source function $f$ is taken as the Ricker wavelet
\begin{equation}
f(t,x,y) = \dfrac{t-2/f_0}{4h^2} \exp\left( -\pi^2 f_0^2 (t-2/f_0)^2 \right) 
\exp\left(\dfrac{(x-0.5)^2+(y-0.5)^2}{4h^2}\right) \text{ for all } (t,x,y) \in \Omega_T,
\end{equation}
where the fine grid parameter and the central frequency are chosen as
$h = 1/256$ and $f_0 = 20$. 
Using the fully discrete scheme, we solve for the numerical solution at the final time $T = 0.2$ 
with time step size $\Delta t = 10^{-4}$. 
We compare the coarse-scale approximation with the fine-grid solution. 
The coarse mesh size varies from $H = 1/64$ to $H = 1/8$, 
and the number of oversampling layers varies according with $m \approx 4 \log(1/H)/ \log(1/8)$. 
In all these combinations, we use $4$ test basis functions per coarse block to construct the 
corresponding localized multiscale basis functions. 

\begin{figure}[ht!]
\centering
\includegraphics[width=0.45\linewidth]{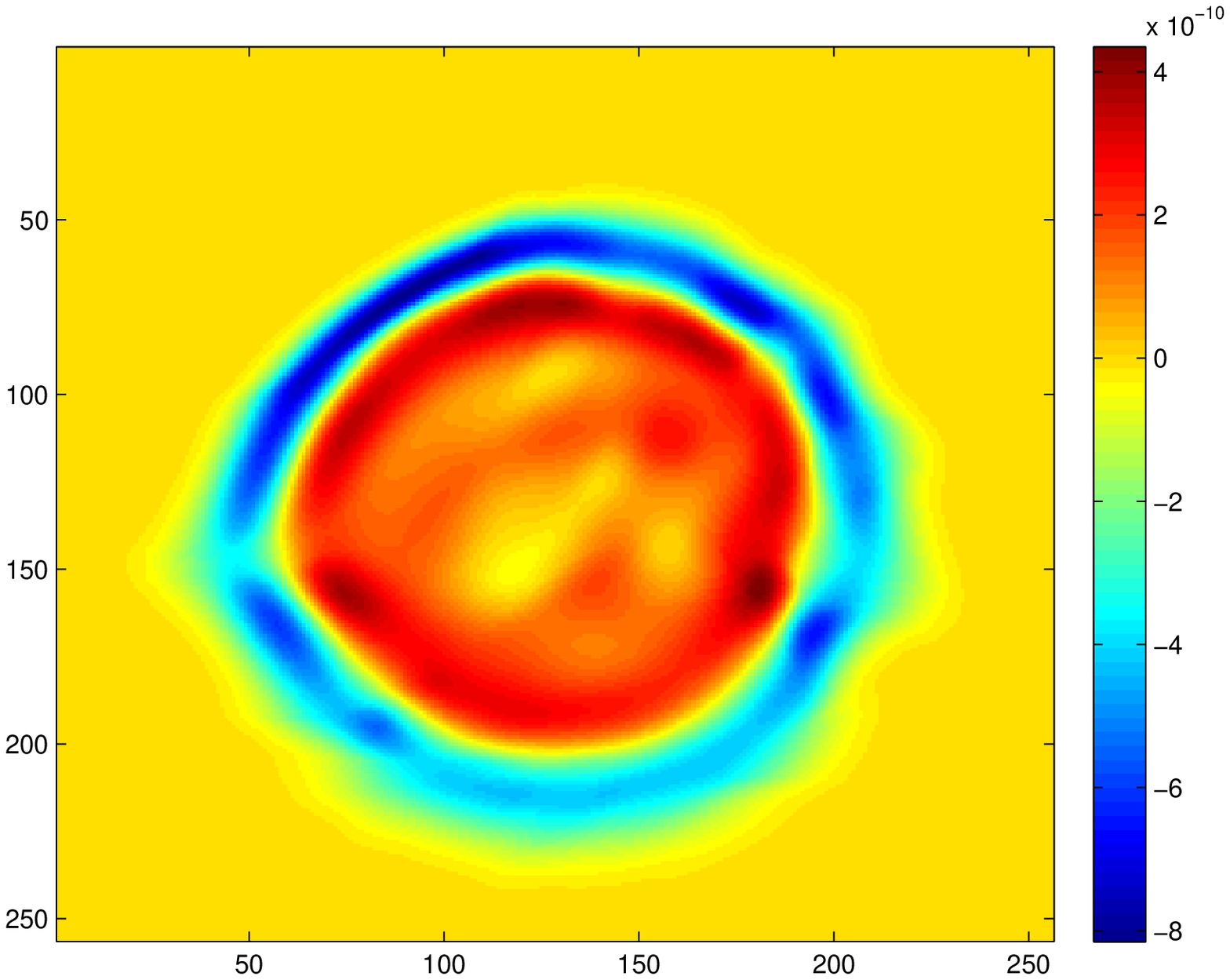}
\includegraphics[width=0.45\linewidth]{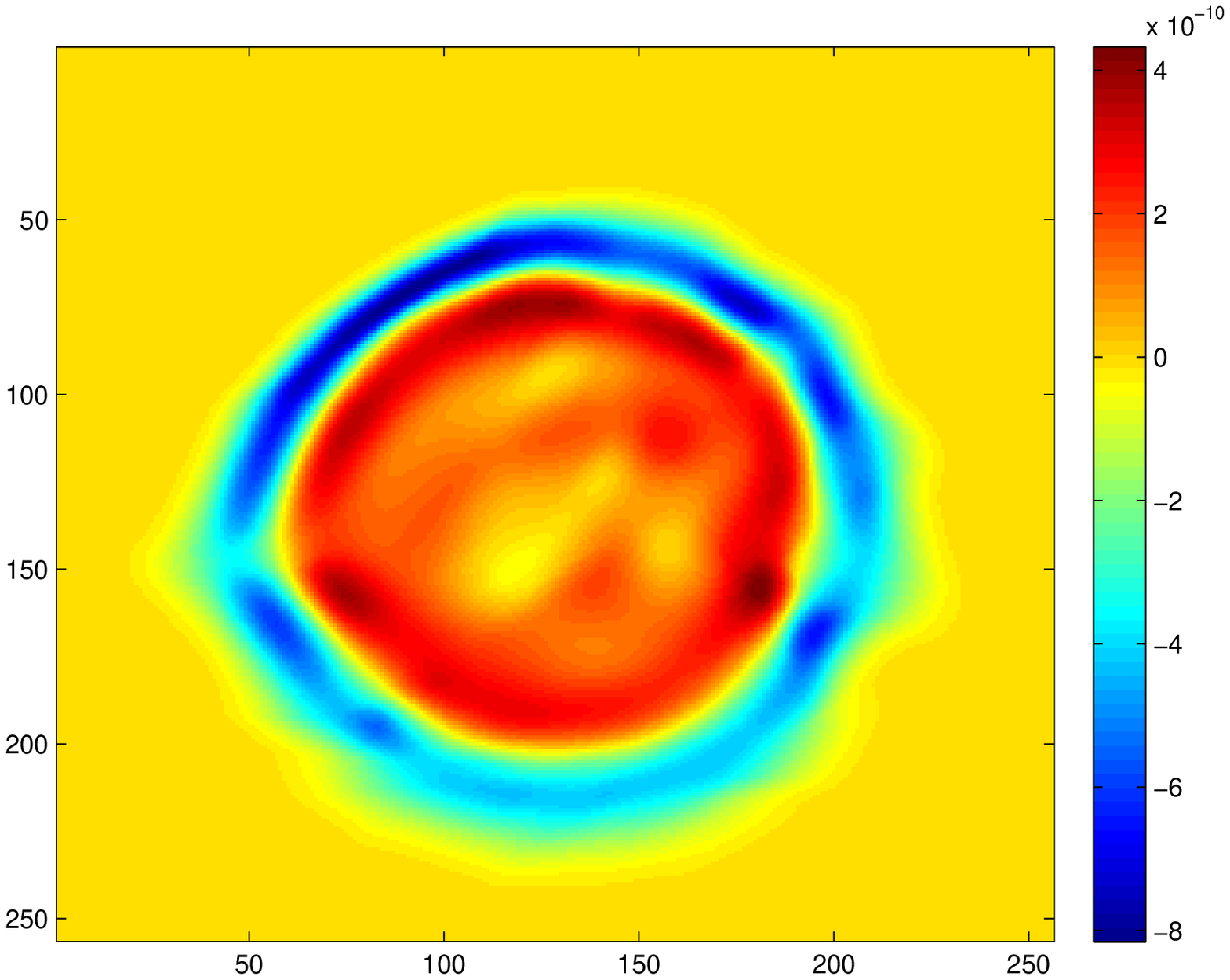}
\caption{Plots of numerical solution for Experiment 2. Fine solution (left) and multiscale solution (right).}
\label{fig:sol_wave}
\end{figure}

Table~\ref{tab:error_wave} records the error of the final solution 
It can been observed that the method results in good 
accuracy and desired convergence in error. 
Figure~\ref{fig:sol_wave} depicts the numerical solutions 
by the fine-scale formulation and the coarse-scale formulation at the final time $T = 0.2$. 
The comparison suggests that our new method provides very good accuracy 
at a reduced computational expense.

\begin{table}[ht!]
\centering
\begin{tabular}{cc|cc}
$m$ & $H$ & Energy error & $L^2$ error \\
\hline
4 & 1/8 & 90.0914\% & 64.3121\% \\
6 & 1/16 & 49.1932\% & 26.4195\% \\
7 & 1/32 & 9.9617\% & 4.4368\% \\
8 & 1/64 & 1.1806\% & 0.5049\% 
\end{tabular}
\caption{History of convergence for wave propagation in Marmousi model.}
\label{tab:error_wave}
\end{table}

\section{Conclusion}
\label{sec:conclusion}
In this paper, we present CEM-GMsDGM, 
a local multiscale model reduction approach in the discontinuous Galerkin framework
for the scalar wave equation. 
The multiscale trial basis functions are defined in coarse oversampled regions by a constraint energy minimization problem, 
which are in general discontinuous on the coarse grid, and coupled by the IPDG formulation 
for solving for a coarse-scale approximation. 
The method is expicit and energy conserving, and 
exhibits both coarse-mesh convergence and spectral convergence, 
provided that the oversampling size is appropriately chosen. 
The stability and the convergence of the method is theoretically analyzed and numerically verified. 
Numerical results for applying the method on a scalar wave equation are also presented.

\section*{Acknowledgements}

The research of Eric Chung is partially supported by the Hong Kong RGC General Research Fund (Project numbers 14304719 and 14302018) and CUHK Faculty of Science Direct Grant 2019-20.

\bibliographystyle{plain}
\bibliography{references}

\begin{thebibliography}{10}

\bibitem{abdulle05}
A.~Abdulle.
\newblock On a priori error analysis of fully discrete heterogeneous multiscale
  fem.
\newblock {\em SIAM J. Multiscale Modeling and Simulation}, 4(2):447--459,
  2005.

\bibitem{abdulle2017multiscale}
Assyr Abdulle and Patrick Henning.
\newblock Multiscale methods for wave problems in heterogeneous media.
\newblock In {\em Handbook of Numerical Analysis}, volume~18, pages 545--576.
  Elsevier, 2017.

\bibitem{arjmand2014analysis}
Doghonay Arjmand and Olof Runborg.
\newblock Analysis of heterogeneous multiscale methods for long time wave
  propagation problems.
\newblock {\em Multiscale Modeling \& Simulation}, 12(3):1135--1166, 2014.

\bibitem{buffa2006analysis}
A~Buffa, TJR Hughes, and G~Sangalli.
\newblock Analysis of a multiscale discontinuous {Galerkin} method for
  convection-diffusion problems.
\newblock {\em SIAM Journal on Numerical Analysis}, 44(4):1420--1440, 2006.

\bibitem{calo2011note}
Victor Calo, Yalchin Efendiev, and Juan Galvis.
\newblock A note on variational multiscale methods for high-contrast
  heterogeneous porous media flows with rough source terms.
\newblock {\em Advances in water resources}, 34(9):1177--1185, 2011.

\bibitem{ch03}
Z.~Chen and T.~Y. Hou.
\newblock A mixed multiscale finite element method for elliptic problems with
  oscillating coefficients.
\newblock {\em Mathematics of Computation}, 72(242):541--576, 2003.

\bibitem{cheung2018mc}
Siu~Wun Cheung, Eric~T Chung, Yalchin Efendiev, Wing~Tat Leung, and Maria
  Vasilyeva.
\newblock Constraint energy minimizing generalized multiscale finite element
  method for dual continuum model.
\newblock {\em arXiv preprint arXiv:1807.10955}, 2018.

\bibitem{cheung2019dg}
Siu~Wun Cheung, Eric~T Chung, and Wing~Tat Leung.
\newblock Constraint energy minimizing generalized multiscale discontinuous
  galerkin method.
\newblock {\em arXiv preprint arXiv:1909.12461}, 2019.

\bibitem{cheung2019bayes}
Siu~Wun Cheung and Nilabja Guha.
\newblock Dynamic data-driven bayesian gmsfem.
\newblock {\em Journal of Computational and Applied Mathematics}, 353:72 -- 85,
  2019.

\bibitem{cgh09}
C.C. Chu, I.G. Graham, and T.~Hou.
\newblock A new multiscale finite element methods for high-contrast elliptic
  interface problem.
\newblock {\em Mathematics of Computation}, 79:1915--1955, 2010.

\bibitem{chung2016adaptive}
Eric Chung, Yalchin Efendiev, and Thomas~Y Hou.
\newblock Adaptive multiscale model reduction with generalized multiscale
  finite element methods.
\newblock {\em Journal of Computational Physics}, 320:69--95, 2016.

\bibitem{chung2018mixed}
Eric Chung, Yalchin Efendiev, and Wing~Tat Leung.
\newblock Constraint energy minimizing generalized multiscale finite element
  method in the mixed formulation.
\newblock {\em Computational Geosciences}, 22(3):677--693, 2018.

\bibitem{chung2017residualWave}
Eric~T Chung, Yalchin Efendiev, Richard~L Gibson, and Wing~Tat Leung.
\newblock Residual-driven online multiscale methods for acoustic-wave
  propagation in 2d heterogeneous media.
\newblock {\em Geophysics}, 82(2):T69--T77, 2017.

\bibitem{chung2015generalizedwave}
Eric~T Chung, Yalchin Efendiev, Richard~L Gibson~Jr, and Maria Vasilyeva.
\newblock A generalized multiscale finite element method for elastic wave
  propagation in fractured media.
\newblock {\em GEM-International Journal on Geomathematics}, pages 1--20, 2015.

\bibitem{chung2014generalized}
Eric~T Chung, Yalchin Efendiev, and Wing~Tat Leung.
\newblock Generalized multiscale finite element methods for wave propagation in
  heterogeneous media.
\newblock {\em Multiscale Modeling \& Simulation}, 12(4):1691--1721, 2014.

\bibitem{chung2015residual}
Eric~T Chung, Yalchin Efendiev, and Wing~Tat Leung.
\newblock Residual-driven online generalized multiscale finite element methods.
\newblock {\em Journal of Computational Physics}, 302:176--190, 2015.

\bibitem{chung2017dg}
Eric~T. Chung, Yalchin Efendiev, and Wing~Tat Leung.
\newblock An online generalized multiscale discontinuous galerkin method
  (gmsdgm) for flows in heterogeneous media.
\newblock {\em Communications in Computational Physics}, 21(2):401–422, 2017.

\bibitem{chung2018constraint}
Eric~T Chung, Yalchin Efendiev, and Wing~Tat Leung.
\newblock Constraint energy minimizing generalized multiscale finite element
  method.
\newblock {\em Computer Methods in Applied Mechanics and Engineering},
  339:298--319, 2018.

\bibitem{chung2016mixedWave}
Eric~T Chung and Wing~Tat Leung.
\newblock Mixed gmsfem for the simulation of waves in highly heterogeneous
  media.
\newblock {\em Journal of Computational and Applied Mathematics}, 306:69--86,
  2016.

\bibitem{chung2018dg}
E.T. Chung, Y.~Efendiev, and W.T. Leung.
\newblock An adaptive generalized multiscale discontinuous galerkin method for
  high-contrast flow problems.
\newblock {\em SIAM Multiscale Modeling and Simulation}, 16(3):1227--1257,
  2018.

\bibitem{chung2014adaptive}
ET~Chung, Y~Efendiev, and G~Li.
\newblock An adaptive {GMsFEM} for high-contrast flow problems.
\newblock {\em Journal of Computational Physics}, 273:54--76, 2014.

\bibitem{ee03}
W.~E and B.~Engquist.
\newblock Heterogeneous multiscale methods.
\newblock {\em Comm. Math. Sci.}, 1(1):87--132, 2003.

\bibitem{emz05}
W.~E, P.~Ming, and P.~Zhang.
\newblock Analysis of the heterogeneous multiscale method for elliptic
  homogenization problems.
\newblock {\em J. Amer. Math. Soc.}, 18(1):121--156, 2005.

\bibitem{egh12}
Y.~Efendiev, J.~Galvis, and T.~Hou.
\newblock Generalized multiscale finite element methods.
\newblock {\em Journal of Computational Physics}, 251:116--135, 2013.

\bibitem{eglmsMSDG}
Y~Efendiev, J~Galvis, R~Lazarov, M~Moon, and M~Sarkis.
\newblock Generalized multiscale finite element method. {Symmetric} interior
  penalty coupling.
\newblock {\em Journal of Computational Physics}, 255:1--15, 2013.

\bibitem{eglp13oversampling}
Y~Efendiev, J~Galvis, G~Li, and M~Presho.
\newblock Generalized multiscale finite element methods. {Oversampling}
  strategies.
\newblock {\em International Journal for Multiscale Computational Engineering,
  accepted}, 2013.

\bibitem{egw10}
Y.~Efendiev, J.~Galvis, and X.H. Wu.
\newblock Multiscale finite element methods for high-contrast problems using
  local spectral basis functions.
\newblock {\em Journal of Computational Physics}, 230:937--955, 2011.

\bibitem{eh09}
Y.~Efendiev and T.~Hou.
\newblock {\em Multiscale Finite Element Methods: Theory and Applications}.
\newblock Springer, 2009.

\bibitem{ehw99}
Y.~Efendiev, T.~Hou, and X.H. Wu.
\newblock Convergence of a nonconforming multiscale finite element method.
\newblock {\em SIAM J. Numer. Anal.}, 37:888--910, 2000.

\bibitem{efendiev2015spectral}
Yalchin Efendiev, Raytcho Lazarov, Minam Moon, and Ke~Shi.
\newblock A spectral multiscale hybridizable discontinuous {G}alerkin method
  for second order elliptic problems.
\newblock {\em Computer Methods in Applied Mechanics and Engineering},
  292:243--256, 2015.

\bibitem{efendiev2017bayes}
Yalchin Efendiev, Wing~Tat Leung, S.~W. Cheung, N.~Guha, V.~H. Hoang, and
  B.~Mallick.
\newblock Bayesian multiscale finite element methods. modeling missing subgrid
  information probabilistically.
\newblock {\em International Journal for Multiscale Computational Engineering},
  15(2):175--197, 2017.

\bibitem{elfverson2013dg}
D.~Elfverson, E.~Georgoulis, A.~Målqvist, and D.~Peterseim.
\newblock Convergence of a discontinuous galerkin multiscale method.
\newblock {\em SIAM Journal on Numerical Analysis}, 51(6):3351--3372, 2013.

\bibitem{engquist2012multiscale}
Bj{\"o}rn Engquist, Henrik Holst, and Olof Runborg.
\newblock Multiscale methods for wave propagation in heterogeneous media over
  long time.
\newblock In {\em Numerical analysis of multiscale computations}, pages
  167--186. Springer, 2012.

\bibitem{hw97}
T.~Hou and X.H. Wu.
\newblock A multiscale finite element method for elliptic problems in composite
  materials and porous media.
\newblock {\em J. Comput. Phys.}, 134:169--189, 1997.

\bibitem{hou2017sparse}
Thomas~Y Hou and Pengchuan Zhang.
\newblock Sparse operator compression of higher-order elliptic operators with
  rough coefficients.
\newblock {\em Research in the Mathematical Sciences}, 4(1):24, 2017.

\bibitem{hfmq98}
T.J.R. Hughes, G.R. Feij\'oo, L.~Mazzei, and J.-B. Quincy.
\newblock The variational multiscale method - a paradigm for computational
  mechanics.
\newblock {\em Comput. Methods Appl. Mech Engrg.}, 127:3--24, 1998.

\bibitem{hughes2007variational}
TJR Hughes and G~Sangalli.
\newblock Variational multiscale analysis: the fine-scale {Green's} function,
  projection, optimization, localization, and stabilized methods.
\newblock {\em SIAM Journal on Numerical Analysis}, 45(2):539--557, 2007.

\bibitem{Iliev_MMS_11}
O.~Iliev, R.~Lazarov, and J.~Willems.
\newblock Variational multiscale finite element method for flows in highly
  porous media.
\newblock {\em Multiscale Model. Simul.}, 9(4):1350--1372, 2011.

\bibitem{jiang2012priori}
Lijian Jiang and Yalchin Efendiev.
\newblock A priori estimates for two multiscale finite element methods using
  multiple global fields to wave equations.
\newblock {\em Numerical Methods for Partial Differential Equations},
  28(6):1869--1892, 2012.

\bibitem{maier2019explicit}
Roland Maier and Daniel Peterseim.
\newblock Explicit computational wave propagation in micro-heterogeneous media.
\newblock {\em BIT Numerical Mathematics}, 59(2):443--462, 2019.

\bibitem{maalqvist2014localization}
Axel M{\aa}lqvist and Daniel Peterseim.
\newblock Localization of elliptic multiscale problems.
\newblock {\em Mathematics of Computation}, 83(290):2583--2603, 2014.

\bibitem{owhadi2017multigrid}
Houman Owhadi.
\newblock Multigrid with rough coefficients and multiresolution operator
  decomposition from hierarchical information games.
\newblock {\em SIAM Review}, 59(1):99--149, 2017.

\bibitem{owhadi2008numerical}
Houman Owhadi and Lei Zhang.
\newblock Numerical homogenization of the acoustic wave equations with a
  continuum of scales.
\newblock {\em Computer Methods in Applied Mechanics and Engineering},
  198(3-4):397--406, 2008.

\bibitem{owhadi2014polyharmonic}
Houman Owhadi, Lei Zhang, and Leonid Berlyand.
\newblock Polyharmonic homogenization, rough polyharmonic splines and sparse
  super-localization.
\newblock {\em ESAIM: Mathematical Modelling and Numerical Analysis},
  48(2):517--552, 2014.

\bibitem{papanicolau1978asymptotic}
G~Papanicolau, A~Bensoussan, and J-L Lions.
\newblock {\em Asymptotic analysis for periodic structures}.
\newblock Elsevier, 1978.

\bibitem{park2019mc}
Jun Sur~Richard Park, Siu~Wun Cheung, Tina Mai, and Viet~Ha Hoang.
\newblock Multiscale simulations for upscaled multi-continuum flows.
\newblock {\em arXiv preprint arXiv:1909.04722}, 2019.

\bibitem{riviere2008discontinuous}
B{\'e}atrice Rivi{\`e}re.
\newblock {\em Discontinuous Galerkin methods for solving elliptic and
  parabolic equations: theory and implementation}.
\newblock Society for Industrial and Applied Mathematics, 2008.

\bibitem{maria2019nlmc}
Maria Vasilyeva, Eric~T Chung, Siu~Wun Cheung, Yating Wang, and Georgy
  Prokopev.
\newblock Nonlocal multicontinua upscaling for multicontinua flow problems in
  fractured porous media.
\newblock {\em Journal of Computational \& Applied Mathematics}, 355:258--267,
  2019.

\bibitem{vdovina2005operator}
Tetyana Vdovina, Susan~E Minkoff, and Oksana Korostyshevskaya.
\newblock Operator upscaling for the acoustic wave equation.
\newblock {\em Multiscale Modeling \& Simulation}, 4(4):1305--1338, 2005.

\bibitem{wang2020vug}
Min Wang, Siu~Wun Cheung, Eric~T. Chung, Maria Vasilyeva, and Yuhe Wang.
\newblock Generalized multiscale multicontinuum model for fractured vuggy
  carbonate reservoirs.
\newblock {\em Journal of Computational and Applied Mathematics}, 366:112370,
  2020.

\bibitem{weh02}
X.H. Wu, Y.~Efendiev, and T.Y. Hou.
\newblock Analysis of upscaling absolute permeability.
\newblock {\em Discrete and Continuous Dynamical Systems, Series B.},
  2:158--204, 2002.

\end{thebibliography}

\end{document}